\documentclass[times,sort&compress,3p]{elsarticle}
\journal{Journal of Multivariate Analysis}
\usepackage[labelfont=bf]{caption}

\usepackage{amsmath,amsfonts,amssymb,amsthm,booktabs,color,epsfig,graphicx,hyperref,url}

\usepackage[linesnumbered,ruled]{algorithm2e}
\usepackage{xparse}
\usepackage{longtable}
\usepackage{array}
\usepackage{multirow}
\usepackage{wrapfig}
\usepackage{colortbl}
\usepackage{pdflscape}
\usepackage{tabu}
\usepackage{threeparttable}
\usepackage{threeparttablex}
\usepackage[normalem]{ulem}
\usepackage{makecell}
\usepackage{bm}
\usepackage{mathtools, bm, bbold}
\usepackage{scalerel}
\usepackage[utf8]{inputenc}
\usepackage{titlesec}
\usepackage{float}
\usepackage{subfig}
\usepackage[dvipsnames]{xcolor}
\theoremstyle{plain}

\newtheorem{corollary}{Corollary}

\newtheorem{property}{Proposition}

\theoremstyle{definition}
\newtheorem{definition}{Definition}

\newtheorem{remark}{Remark}
\newtheorem{hypothesis}{Hypothesis}
\newtheorem{dataset}{Dataset}

\DeclareMathOperator*{\argmin}{{\rm arg\,min}}
\DeclareMathOperator*{\Times}{\times}

\makeatletter
\newcommand{\@giventhatstar}[2]{\ensuremath{\left(#1\;\middle|\;#2\right)}}
\newcommand{\@giventhatnostar}[3][]{#1(#2\;#1|\;#3#1)}
\newcommand{\giventhat}{\@ifstar\@giventhatstar\@giventhatnostar}
\makeatother

\renewcommand{\L}{\mathcal{L}}
\newcommand{\norme}[1]{\lVert #1\rVert_{2}^{2}}
\newcommand{\scalarprod}[2]{\,\langle #1,#2\rangle}
\newcommand{\normeL}[1]{\lVert #1\rVert_{\L}^{2}}
\newcommand{\scalarprodL}[2]{\,\langle #1,#2\rangle_{\L}}

\newcommand{\indep}{\perp\!\!\!\perp}
\newcommand{\distL}[2]{d_{\L}(#1,#2)^2}

\newcommand{\Proba}[1]{\Pr\left(#1\right)}
\newcommand{\Esper}[1]{\mathbb{E}\left(#1\right)}

\newcommand{\Gras}[1]{\bm{#1}}

\newcommand{\Emph}[1]{#1} % should be \emph{#1}
\newcommand{\Bold}[1]{#1} % should be \textbf{#1}

\renewenvironment{proof}{{\bfseries Proof}}{\qed}

\parskip = \baselineskip

\begin{document}

\begin{frontmatter}

\title{Dependence structure estimation using Copula Recursive Trees}

\author[A1]{Oskar Laverny\corref{mycorrespondingauthor}}
\author[A1]{Esterina Masiello}
\author[A1]{Véronique Maume-Deschamps}
\author[A2]{Didier Rullière}
%\affil[1]{Université Claude Bernard Lyon 1, France}
%\affil[2]{SCOR SE, France}
%\affil[3]{Insitut Camille Jordan UMR 5208, France}
%\affil[4]{Laboratoire SAF EA 2429, France}
%
%
%\author[A1]{Author One (First name, then family name)}
%\author[A2]{Author Two\corref{mycorrespondingauthor}}

\address[A1]{Institut Camille Jordan, UMR 5208, Université Claude Bernard Lyon 1, Lyon, France}
\address[A2]{Mines Saint-Etienne, Univ Clermont Auvergne, CNRS, UMR 6158 LIMOS, Institut Henri Fayol, F - 42023 Saint-Etienne France}
\cortext[mycorrespondingauthor]{Corresponding author. Email address: \url{oskar.laverny@univ-lyon1.fr}}

\begin{abstract}
We construct the COpula Recursive Tree (CORT) estimator:  a flexible, consistent, piecewise linear estimator of a copula, leveraging the patchwork copula formalization and various piecewise constant density estimators. While the patchwork structure imposes a grid, the CORT estimator is data-driven and constructs the (possibly irregular) grid recursively from the data, minimizing a chosen distance on the copula space. The addition of the copula constraints makes usual density estimators unusable, whereas the CORT estimator is only concerned with dependence and guarantees the uniformity of margins. Refinements such as localized dimension reduction and bagging are developed, analyzed, and tested through simulated data.
\end{abstract}


\end{frontmatter}

\section{Introduction\label{sec:1}}

Although the estimation of copula~\cite{nelsen2006,joe2014,durante2015a} is a wide-treated subject, most efficient estimators available in the literature are based on restricted, parametric estimation. Vine copulas~\cite{nagler2016,nagler2018,muller2017,muller2019,okhrin2017}, although useful in high dimensions, often use parametric models, such as Archimedean copulas, as base building blocks. On the other hand, graphical models~\cite{li2018,friedman2008} assume a Gaussian dependence structure and therefore are fast but under restrictive assumptions. Classical nonparametric density estimators such as kernels~\cite{hickernell1998,siloko2016,siloko2018,siloko2019} or  wavelets~\cite{genest2012,morettin2010,gavish2010} are not suited to satisfy constraints such as the uniformity of margins (one counter-example may be found in~\cite{chen2007,fermanian2004}). We explore here a specific class of non-parametric copula density estimators with tree-structured piecewise constant densities, and design an estimator that lies in this class, the CORT estimator.

The CORT estimator is based on the \Emph{density estimation tree} from~\cite{ram2011}, a tree-structured non-parametric density estimator, and on the framework of \Emph{patchwork copulas} from~\cite{durante2015,deamo2017,durante2013}. There already exist several other piecewise constant density estimators: the \Emph{cascaded histograms} of~\cite{goh2015}, the \Emph{Dirichlet-based Polya tree}~\cite{ning2018}, the \Emph{distribution element trees} by~\cite{meyer2018a}, the \Emph{adaptative sparse grids} of~\cite{peherstorfer2014}, the framework of \Emph{Gaussian mixtures} by~\cite{criminisi2012}, the \Emph{Bayesian sequential partitioning} techniques by~\cite{lin2016,li2018a} with their interesting asymptotic consistency results, and the \Emph{Wasserstein compression} techniques provided by~\cite{luini2020} are all worth noting in the field of non-parametric piecewise density estimation. But these models are built to estimate densities without taking into account uniformity of margins, and they do not always lead to proper copulas when applied on pseudo-observations or true copula samples.

The CORT estimator has the particularity of being tree-shaped which ensures on one hand that the computation of the estimated density and the distribution function on new data points is fast, and on the other hand that the storage of the model is efficient. Thus, it could be used for tasks such as re-sampling a dataset outside the already existing points, or for compression purposes, when dealing with big-data dependencies. Finally, under mild conditions, the estimator is a proper copula, where classical non-parametric estimators, such as Deheuvel's empirical copula, are not.

This paper is organized as follows. Section~\ref{sec:2} describes the class of piecewise linear copulas and gives some of their properties. In Section~\ref{sec:3}, we propose an estimation procedure, allowing localized dimension reduction, and we establish a convergence result for this procedure. Section~\ref{sec:4} deals with ensemble models based on the CORT estimator: bagging techniques and out-of-bag generalization statistics are developed in the field of copula density estimation, and applied to the CORT estimator. Finally, Section~\ref{sec:5} investigates the performance of the model by applications on some simulated examples, and Section~\ref{sec:6} concludes.

\section{The piecewise linear copula\label{sec:2}}

Let $\Gras{X} = (X_1,\ldots,X_d)^{\top}$ be a multivariate random vector of dimension $d$. We are interested in the dependence structure of $\Gras{X}$. The concept of copula, whose formalization is due to~\cite{sklar1959}, allows to study this dependence separately from the marginal distributions. Consider the distribution function (d.f.) $F$ of the random vector $\Gras{X}$ with marginal d.f.s $F_1,\ldots,F_d$, and define the copula $C$ as: $$C(\Gras{u}) = F\left\{F_1^{-1}(u_1),\ldots,F_d^{-1}(u_d)\right\},\quad \forall \Gras{u} \in [0,1]^d,$$

where $F_{i}^{-1}: y \mapsto \inf\;\left\{x \in \bar{\mathbb{R}}: F_{i}(x) \ge y\right\}$ is the generalized inverse of $F_{i}$. Then the distribution $F$ of $\Gras{X}$ is characterized by $F_1,\ldots,F_d$ and $C$. Sklar's Theorem~\cite{sklar1959} states that a copula $C$ satisfying the previous equality always exists, and that it is uniquely defined provided that the marginal random variables are absolutely continuous. In particular, it is unique if the random vector is continuous. Note that $C$ is the distribution function of a $d$-dimensional random vector with uniform margins.

The estimation of the full distribution can then be split into the estimation of one-dimensional margins, a widely treated subject, and the estimation of the copula for which we propose here to work with a piecewise linear distribution function. In the following, we define the piecewise linear copulas and then present some of their properties.
	 
	\subsection{Definition}
	
	Let $\mathbb{I} = [0,1]^d$ be the unit hypercube, the domain of any copula, and let $\mathbb 1_{A}$ be the indicator function of an event $A$.
	
	\begin{definition}[Piecewise linear copula]
		
		Let $\L$ be a finite partition of $\mathbb{I}$ into subsets called \Emph{leaves}, usually denoted by $\ell$. The piecewise linear copula with partition $\L$ and weighs $p$ is defined by its distribution function: 
		
		\begin{equation}\label{eq:definition_copula}
		\forall \Gras{u} \in \mathbb{I},\; C_{\Gras{p},\L}(\Gras{u}) = \sum\limits_{\ell \in \L} p_\ell \lambda_\ell(\Gras{u}),
		\end{equation}
		
		where $\lambda$ denotes the Lebesgue measure of a set, $\lambda_\ell(\Gras{u}) = \lambda(\ell)^{-1}\lambda([0,\Gras{u}]\cap \ell)$, and $\Gras{p}$ is a vector of non-negative weights summing to one. The corresponding copula density, which is piecewise constant, is given by: 
		\begin{equation}\label{eq:patchwork_density}
		c_{\Gras{p},\L}(\Gras{u}) = \sum\limits_{\ell \in \L}  \frac{p_\ell}{\lambda(\ell)}\mathbb{1}_{\Gras{u} \in \ell}.
		\end{equation}
		
	\end{definition}
	
	This type of histogram, where leaves might not all have same shape and size, has already been used in the case of density estimation with different construction schemes and leaves shapes, e.g. with a Voronoï diagram~\cite{cervellera2020,geiss2013} or a Delaunay tessellation~\cite{watson1981,hornus1992,boissonnat2019} as partition, or more trivially with simple sets of hyper-boxes~\cite{ram2011,anderlini2016,lin2016,luini2020}.
	
	\begin{remark}[Existence]\label{rem:existance} Depending on the choice of the partition $\L$ and the weights $\Gras{p}$, the distribution function $C_{\Gras{p},\L}$ is not always a copula. However, if $\forall \ell \in \L,\; p_\ell = \lambda(\ell)$, then $C_{\Gras{p},\L}$ is the independence copula. Therefore, for any partition, there exists at least one set of weights making the model a proper copula. On the other hand, for a partition that is too complex, this is frequently the only solution: the assumption of uniform leaves and the marginals uniformity constraints restrict the shape of bins. For polytopic leaves that are not hyper-rectangles, we do not know if weighting them efficiently can be achieved inside the copula constraints.
	\end{remark}
	Before looking more precisely at the copula constraints on these piecewise constant distribution functions, we therefore restrict ourselves to the case of hyper-rectangular leaves, leading to the following definition:
	
	\begin{definition}[Hyper-rectangles and suitable partitions] Let $\Gras{a}$ and $\Gras{b}$ both be in $\mathbb{I}$. Then, if $\Gras{a} \le \Gras{b}$ (component-wise), we define the hyper-rectangle $(\Gras{a},\Gras{b}]$ as $(\Gras{a},\Gras{b}] = (a_1,b_1] \times ... \times (a_d,b_d]$. We call \Emph{suitable} a partition where every leave is a hyper-rectangle with strictly positive Lebesgue measure.
	\end{definition}
	Remark~\ref{rem:existance} partly drives the definition of a suitable partition. It is also why we chose to extend the \Emph{density estimation trees} from~\cite{ram2011} instead of another piecewise constant density estimator: it produces a suitable partition. 
	If not specified, we consider by default that partitions we are dealing with are suitable. In the next subsection, we propose some properties of the dependence structure induced by such a copula. 
	
	\subsection{Properties}
	
	With the above formulation of a piecewise linear copula, we can easily obtain closed-form expressions for classical quantities of interest in copula modeling. We recall some of those quantities and then derive their expression for piecewise linear copulas. The Kendall $\tau$ and Spearman $\rho$ (see~\cite{nelsen2006}) are common dependence measures that can be computed from a copula. They are respectively defined for a copula $C$ and its density $c$ as:
	$$\tau = 4 \int C(\Gras{u}) \, c(\Gras{u}) \;d\Gras{u} -1,\quad\rho = 12 \int C(\Gras{u}) d\Gras{u} -3.$$
	Note that both $\tau$ and $\rho$ are always in $[-1,1]$. The piecewise constant expression of the density  in the piecewise linear class allows for simple computation of $\tau$ and $\rho$, although the expressions can be somewhat cumbersome.
	
	\begin{property}[Common dependence measures]\label{ppt:comon_dep_mes}
		Let $C_{\Gras{p},\L}$ be a piecewise linear copula. Its Kendall $\tau$ and Spearman $\rho$ are given in closed form by: \begin{align*}
		\tau &= -1 + 4\sum\limits_{\substack{\ell \in \L\\ k \in \L}}\frac{p_\ell p_k}{\lambda(\ell) \lambda(k)}\prod\limits_{i=1}^{d} \left\{\frac{(b\wedge d)^2 - (a\vee c)^2}{2} + c\left((a\vee c) - (b\wedge d)\right)+(d-c)\left(b - (a\vee d)\right)\right\},\\
		\rho &= -3 + 6\sum\limits_{\ell \in \L} p_\ell \prod\limits_{i=1}^{d} \left(2 - b_i - a_i\right),
		\end{align*}
		
		where we denote $\ell = (\Gras{a},\Gras{b}]$ and $k = (\Gras{c},\Gras{d}]$, and $\wedge,\;\vee$ denote respectively the minimum and maximum operator.		
	\end{property}
	The proof is postponed to the appendix. Matrices of bivariate dependence measures can be obtained by projection of the partition on couples of dimensions and using the same formula on the projected models.
	The closed form expression for piecewise linear copulas facilitates greatly these kinds of computations. This property will be exploited later, when introducing penalization techniques. Furthermore, these closed form expressions will be used to assess the performance of the fitting procedure that we describe in the next section.

\section{Estimation\label{sec:3}}

Suppose that we have a dataset $\left(u_{i,j}\right)_{n\times d}$ of (pseudo-)observations from an unknown copula $C$. We seek parameters $(\Gras{p},\L)$ of $c_{\Gras{p},\L}^{(n)}$, an approximation of $c$ in the piecewise linear copula class based on these $n$ observations. To find the optimal parameters $(\Gras{p}^{*},\L^{*})$, we will adopt a two stage mechanism, considering first that the partition $\L$ is known. From now on, we denote by $\norme{f}$ the squared $L_2$ norm of a function $f$, given by $\norme{f} = \int f(x)^2 dx$.

\subsection{Optimal weights \texorpdfstring{$\Gras{p}^{*}_{\L}$}{pL*} knowing the partition \texorpdfstring{$\L$}{L}}

Suppose that a partition $\L$ is already constructed, and that we want to construct weights $\Gras{p}$ to complete the approximation. As was done by~\cite{ram2011} (see also~\cite{alquier2008,birge2013}), we will use an Integrated Square Error (ISE) loss to build the weights. Given a copula density $c$, the ISE of an estimator $\hat{c}$ is the squared $L_2$ distance to $c$: $$I(\hat{c}) = \norme{\hat{c} - c} = \int (\hat{c}(\Gras{u}) - c(\Gras{u}))^2 d\Gras{u}.$$

To approximate the true copula density $c$ by a piecewise linear copula, we are looking to solve:
\begin{align}
\argmin\limits_{\Gras{p},\L}\; I(c_{\Gras{p},\L}^{(n)}) &= \argmin\limits_{\Gras{p},\L}\;\norme{c_{\Gras{p},\L}^{(n)} - c} \nonumber =\argmin\limits_{\Gras{p},\L}\;\int (c_{\Gras{p},\L}^{(n)}(\Gras{u}) - c(\Gras{u}))^2 d\Gras{u} \nonumber\\
&=\argmin\limits_{\Gras{p},\L}\;\int c_{\Gras{p},\L}^{(n)}(\Gras{u})^2 - 2 c(\Gras{u})c_{\Gras{p},\L}^{(n)}(\Gras{u}) d\Gras{u} \nonumber = \argmin\limits_{\Gras{p},\L}\; \norme{c_{\Gras{p},\L}^{(n)}} - 2 \scalarprod{c_{\Gras{p},\L}^{(n)}}{c},\label{eq:ise_argmin}
\end{align}

since $\norme{c}$ is irrelevant to the minimization. Remark that $\scalarprod{c_{\Gras{p},\L}^{(n)}}{c} = \Esper{c_{\Gras{p},\L}^{(n)}(U)}$, so that, with a slight abuse of notation, we write $\scalarprod{c_{\Gras{p},\L}^{(n)}}{c}$ even if $C$ does not admit a density $c$.
Furthermore, since the true copula $C$ and its density $c$ are unknowns, using empirical observations $(u_{i,j})$ from the copula, $\Esper{c_{\Gras{p},\L}^{(n)}(U)}$ can be approximated by $n^{-1}\sum_{i=1}^n c_{\Gras{p},\L}^{(n)}(u_i)$. We then define our empirical loss.

\begin{definition}[Empirical Integrated Square Error] Given observations $(u_{i,j})$ from a copula, define the Empirical Integrated Square Error (EISE) of an estimator $\hat{c}$ of the copula density as: \begin{equation}\label{eq:def_eise}
\hat{I}(\hat{c}) = \norme{\hat{c}} - \frac{2}{n}\sum_{i=1}^n \hat{c}(\Gras{u_i}).
\end{equation}
	
\end{definition}

To minimize this loss, we find the weights $\Gras{p}_{\L}^{*}$ knowing the partition $\L$ by solving the following problem: \begin{align}\label{eq:ise_argmin_empirical} \argmin\limits_{p}\hspace{0.5em} &\hat{I}(c_{\Gras{p},\L}),\quad\mathrm{s.t.}\hspace{0.5em} c_{\Gras{p},\L}\text{ is a copula.}
\end{align}

The copula constraints in Eq. \eqref{eq:ise_argmin_empirical} are classically expressed in the literature as in Definition~\ref{def:raw_constraints}:

\begin{definition}[Copula constraints]\label{def:raw_constraints} Let $C$ be the distribution function of a signed measure $\mu_C$. $C$ is a copula if the following three conditions are satisfied:\begin{align*}
	&\forall \Gras{u} \in \mathbb{I}, \forall i \in \{1,\ldots,d\}, C(u_1,\ldots,u_{i-1},0,u_{i+1},\ldots,u_d) = 0, &&\text{(Ground constraint)},\\
	&\forall \Gras{u} \in \mathbb{I}, \forall i \in \{1,\ldots,d\}, C(1,\ldots,1,u_i,1,\ldots,1) = u_i,&&\text{(Marginal uniformity)},\\ 
	&\forall \Gras{a} \le \Gras{b} \in \mathbb{I}, \mu_{C}\left([\Gras{a},\Gras{b}]\right) \ge 0,  &&\text{($d$-increasingness)}.
	\end{align*}
\end{definition}

Note that, if $C$ is a distribution function of some random vector, the first and third conditions are verified. It turns out that, under these constraints, our optimization problem is in fact a quadratic program, as Proposition~\ref{prop:final_quadratic_program} shows.

\begin{property}[Quadratic program]\label{prop:final_quadratic_program}
	Let $\L$ be a suitable partition. Then the weights minimizing the empirical integrated square error \eqref{eq:def_eise} under the copula constraints from Definition~\ref{def:raw_constraints} are the unique solution of the quadratic program:
	\begin{align*}
	\argmin\limits_{\Gras{p}}   \hspace{0.5em} &\Gras{p}' \Gras{A}_{\L} \Gras{p}  - 2 \Gras{p}' \Gras{A}_{\L} \Gras{\mathrm{f}}_{\L},\quad\mathrm{s.t.} \hspace{0.5em} \Gras{B}_{\L}\Gras{p}=\Gras{g}_{\L} \text{ and }\Gras{p} \geq \Gras{0},
	\end{align*}

	where the matrices $\Gras{A}_{\L},\;\Gras{B}_{\L}$ and the vectors $\Gras{\mathrm{f}}_{\L},\;\Gras{g}_{\L}$ are given by:
	\begin{align*}
	\Gras{A}_{\L} &= \left(\lambda(\ell)^{-1}\mathbb{1}_{\ell = k}\right)_{\ell \in \L, k \in \L}, && \text{(size } \vert \L \vert\times \vert \L \vert\text{)},\\
	\Gras{B}_1 &= \left(\lambda_{\ell_{i}}(u_i)\right)_{(i,\Gras{u}) \in \{1,\ldots,d\} \times M_{\L}, \; \ell \in \L}, && \text{(size } nd\times \vert \L \vert\text{)},\\
	\Gras{B}_{\L} &= (\Gras{B}_1,\Gras{B}_2), && \text{(size } (nd+1)\times \vert \L \vert\text{)},\\
	\Gras{\mathrm{f}}_{\L} &= \left(\frac{1}{n}\sum\limits_{i=1}^{n} \mathbb{1}_{\Gras{u_i} \in \ell}\right)_{\ell \in \L},&& \text{(size }\vert \L \vert\text{)},\\
	\Gras{B}_2 &= \left(1\right)_{\ell \in \L} ,&& \text{(size } 1\times \vert \L \vert\text{)},\\
	\Gras{g}_1 &= \left(u_i\right)_{(i,\Gras{u}) \in \{1,\ldots,d\} \times M_{\L}} ,&& \text{(size }nd\text{)},\\
	\Gras{g}_{\L} &= (\Gras{g}_1,1) ,&& \text{(size }(nd+1)\text{)},
	\end{align*}
	
	where we denoted $\vert\L\vert$ the cardinal of $\L$, and $M_{\L}$ the set of middle-points of leaves. 
	
	\begin{proof} This proof is in three parts. First, we show that the objective formulation is correct, then we discuss the constraints formulation and finally we prove existence and uniqueness of a solution. Introduce first a new scalar product on $\mathbb{R}^{\vert\L\vert}$: $$\scalarprodL{\Gras{x}}{\Gras{y}} = \sum_{\ell\in\L} \frac{x_\ell y_\ell}{\lambda(\ell)}.$$ We denote by $\normeL{.}$ its associated square norm, and by $d_{\L}$ its associated distance. Note that the associated bilinear symmetric operator has matrix $\Gras{A}_{\L}$ (defined above), a diagonal and positive definite matrix.
		
		Using the definition of $\Gras{A}_{\L}$, $\Gras{\mathrm{f}}_{\L}$ and $\scalarprodL{.}{.}$, the objective function in Eq.~\eqref{eq:ise_argmin_empirical} rewrites:
		
		\begin{equation}
		\hat{I}(c_{\Gras{p},\L}) = \normeL{\Gras{p}} - 2\scalarprodL{\Gras{p}}{\Gras{\mathrm{f}}_{\L}} = \Gras{p}' \Gras{A}_{\L} \Gras{p}  - 2 \Gras{p}' \Gras{A}_{\L}\Gras{\mathrm{f}}_{\L}.\label{eq:expression_prod_scalar_of_the_unconstrainted_loss}
		\end{equation}
		
		Then, the unconstrained version of the convex optimization problem from Eq.~\eqref{eq:ise_argmin_empirical} corresponds to the projection of $\Gras{\mathrm{f}}_{\L}$ onto $[0,1]^{\vert\L\vert}$ (convex, closed), with respect to the norm $\normeL{.}$. Note that, if we set the weights to be equal to the empirical frequencies $\Gras{\mathrm{f}}_{\L}$, this result yields, for the optimization of $\L$, a loss of the form: 
		\begin{equation}\label{eq:loss_without_constraint}
		- \normeL{\Gras{\mathrm{f}}_{\L}}  = -\sum\limits_{\ell \in \L}  \frac{\mathrm{f}_\ell^2}{\lambda(\ell)},
		\end{equation}
		
		which is the loss that was directly used by~\cite{ram2011}. 
		
		We now discuss the constraints. Denote by $\mathcal{C}_{\L}$ the subset of $[0,1]^{\vert\L\vert}$ containing vectors $\Gras{p}$ such that $C_{\Gras{p},\L}$ satisfies the set of constraints from Definition~\ref{def:raw_constraints}, for a given $\L$. The claim is that: $$\mathcal{C}_{\L} = \left\{\Gras{p} \in \mathbb{R}^{\vert\L\vert}: \; \Gras{B}_{\L}\Gras{p} = \Gras{g}_{\L} \text{ and }\Gras{p}\ge 0\right\}.$$
		
		The first constraint of Definition~\ref{def:raw_constraints} is trivially satisfied by our model. We show that the second constraint can be evaluated on only one point per leaf. Remember that the piecewise linear copula is uniformly distributed on each leaf. Hence, for all $i \in \{1,\ldots,d\}$, if on some point $\Gras{u} = (u_1,\ldots,u_d)^{\top}$, $\sum_{\ell \in \L} p_{\ell} \lambda_{\ell_{i}}(u_i) = u_i$, then defining $\Gras{x}$ such that $\Gras{u}+\Gras{x}$ is in the same leaf as $\Gras{u}$ will give us that  $\sum_{\ell \in \L} p_{\ell} \lambda_{\ell_{i}}(u_i + x_i) = u_i + x_i$ since only one leaf is active in the sum. Therefore, the marginal uniformity constraint may be evaluated on only one point per leaf, and hence is equivalent to: 
		\begin{equation*}
		\forall \ell \in \L, \exists\; \Gras{u} \in \ell, \forall i \in \{1,\ldots,d\}, \sum\limits_{\ell \in \L} p_\ell \lambda_{\ell_{i}}(u_i) = u_i,
		\end{equation*}
		
		where for $\ell = [\bm a, \bm b]$ we still denote $\ell_i = [a_i, b_i]$ its marginals and $\lambda_{\ell_i}(u_i) = \frac{\lambda([0,u_i]\cap\ell_i)}{\lambda(\ell_i)}$. Then, if we choose evaluation points as middle-points of leaves to put the constraints in matrix-vector form, we have the expression $\Gras{B}_1\Gras{p}=\Gras{g}_1$ for these constraints. Furthermore, we need to force the sum of weights $\Gras{p}$ to be equal to $1$ (so that the total marginal mass is $1$), giving a last line of ones to $\Gras{B}_{\L}$ and a last value of one to $\Gras{g}_{\L}$.

		The third constraint states that the measure associated with the copula $C$ has to be positive on any hyper-rectangle $[\Gras{a},\Gras{b}]$. Recall that:
		$$\forall \Gras{a} \le \Gras{b},\; \mu_C([\Gras{a},\Gras{b}]) = \int\limits_{[\Gras{a},\Gras{b}]} c(\Gras{u})d\Gras{u} = \sum\limits_{\ell \in \mathcal L} p_\ell\lambda_{\ell}([\Gras{a},\Gras{b}])$$ If all weights are positive, then $\forall \Gras{a} \le \Gras{b},\; \mu_C([\bm a,\bm b]) \ge 0$. On the other hand, if one weight is negative, then taking $\Gras{a}, \Gras{b}$ inside the corresponding leaf would make $\mu_C([a,b])$ negative. The last constraint is therefore equivalent to positivity of weights, which gives the last part of the wanted expression for $\mathcal{C}_{\L}$. 
		Finally, existence and uniqueness of a solution are trivial since the objective function is a quadratic function, $\mathcal{C}_{\L}$ is a closed and convex set, and $\mathcal{C}_{\L}$ is non-empty by Remark~\ref{rem:existance}.	\end{proof}
\end{property}

%
%
%\begin{lemma}\label{lem:final_loss} Let $\L$ be a partition of $\mathbb{I}$. Then the weights $\Gras{p}$, solution of \eqref{eq:ise_argmin_empirical} are given as the unique solution of the quadratic program with objective:  \begin{equation}\label{eq:ise}\Gras{p}' \Gras{A}_{\L} \Gras{p}  - 2 \Gras{p}' \Gras{A}_{\L} \Gras{\mathrm{f}}_{\L}\end{equation}
%	
%	
%	 Note that $\Gras{A}_{\L}$ depends only on the partition and that $\Gras{\mathrm{f}}_{\L}$ represents the empirical frequencies in the leaves.
%	 
%	 \begin{proof} 
%		
%	\end{proof}
%\end{lemma}
%
%
%
%\begin{lemma}\label{lem:final_constraints}  The set $\mathcal{C}_{\L}$ is closed, convex, non-empty and writes:
%	
%	
%	
%	
%	\begin{proof}
%		
%	\end{proof}
%\end{lemma}
%Lemmas~\ref{lem:final_loss} and~\ref{lem:final_constraints} lead to the following property, summarizing the previous ideas. 

Denoting $P_{\L,\mathcal{S}}(\Gras{x})$ the orthogonal projection of a vector $\Gras{x}$ onto a set $S$ regarding the distance $d_\L$, the quadratic program from Proposition~\ref{prop:final_quadratic_program} gives the optimal weights knowing the partition as: $$\Gras{p}_{\L}^{*} = P_{\L,\mathcal{C}_{\L}}(\Gras{\mathrm{f}}_{\L}).$$ The empirical frequencies $\Gras{\mathrm{f}}_{\L}$, which are the unconstrained solution, can then be used as a good starting point for a solver. We concentrate now on the construction of the partition $\L$.

\subsection{Locally optimal splitting point \texorpdfstring{$\Gras{x}_{\ell,\Gras{D}}^{*}$}{xlD*}}

Suppose that we have already a suitable partition $\L$ and associated weights $\Gras{p}_{\L}$ such that $C_{\Gras{p}_{\L},\L}$ is a copula. For a given leaf $\ell \in \L$, denote $L = {\ell_1,\ldots,\ell_k}$ a partition of $\ell$ into $k$ new leaves such that $\L_2 = \L\setminus \{\ell\} \cup L$ is a new suitable partition. Then we have, as in~\cite{ram2011}, the following property: 

\begin{property}[Independence of surrogate loss]
	Define the surrogate loss associated to the additional split from $(\Gras{p}_{\L}^{*},\L)$ to $(\Gras{p}_{\L_2}^{*},\L_2)$ as the difference of empirical integrated squared errors $\hat{I}(c_{\Gras{p}_{\L_2}^{*},\L_2}) - \hat{I}(c_{\Gras{p}_{\L}^{*},\L}).$ Then the surrogate loss depends only on the partition $L$ and data inside it.

	\begin{proof} Index the objects of Proposition~\ref{prop:final_quadratic_program} by the partitions $\L$ and $\L_2$, and the part corresponding to $\L \setminus \{\ell\}$ cancels.
	\end{proof}
\end{property}
This locality of the loss allowed~\cite{ram2011} to use a recursive partitioning algorithm. We then only perform \Emph{simple splits}.

\begin{definition}[Simple split and splitting dimensions]
	Denote $\Gras{x}$ a given breakpoint in the leaf $\ell$ and $\Gras{D} \subseteq \{1,\ldots,d\}$ the set of splitting dimensions. Then the simple split of $\ell$ on $\Gras{x}$ with dimensions $\Gras{D}$ is defined as the partition $L(\ell,\Gras{x},\Gras{D})$ given by: 
	$$L((\Gras{a},\Gras{b}],\Gras{x},\Gras{D}) = \Times\limits_{j\in \Gras{D}} \left\{(a_j,x_j],(x_j,b_j]\right\} \Times\limits_{j \in \{1,\ldots,d\}\setminus \Gras{D}} \left\{(a_j,b_j]\right\}\text{.}$$
	
	The full partition of $\mathbb{I}$ obtained after the split is denoted $\L_{\Gras{x},\Gras{D}} = \L \setminus \{\ell\} \cup L(\ell,\Gras{x},\Gras{D})$. When $\Gras{D} = \{1,\ldots,d\}$, we might omit it in the subscripts.
\end{definition}

\begin{remark}[Degrees of freedom]
	In a simple split on a set of dimension $\Gras{D}$, the weighting of the new leaves is a quadratic program with $2^{\vert \Gras{D} \vert}$ parameters responding to $\vert \Gras{D}\vert+1$ linear constraints. Hence, there exists a trade-off between complexity and expressiveness of the model in the dimension $\vert \Gras{D} \vert$ of the breakpoints. We will exploit this characteristic of the recursive procedure for sparsity purposes later on.
\end{remark}

\begin{remark}[No one-dimensional splits]
	Note that the copula constraints will not allow for estimation with only one-dimensional splits ($\vert \Gras{D} \vert = 1$), as in~\cite{ram2011}, since there would be no degrees of freedom in the weights. This represents a huge problem as multivariate splits often imply bigger computational burden. Furthermore, the constraints themselves are \Emph{not} localized, but on the global scale, hence including them forbids a parallel implementation. We will see later that this issue can be avoided by delaying the constraint problem to a later stage of the optimization process. 
\end{remark}

Note that, neglecting the constraints, knowing $\Gras{D}$, we are going to choose the splitting point $\Gras{x}_{\ell,\Gras{D}}^{*}$ as:

\begin{equation}\label{eq:surrogate_loss_last_version}
\Gras{x}_{\ell,\Gras{D}}^{*} = \argmin\limits_{\Gras{x}\in\ell}  - \sum\limits_{k \in L(\ell,\Gras{x},\Gras{D})} \frac{\mathrm{f}_k^2}{\lambda(k)}
\end{equation}

In the next subsection, before giving a complete description of the fitting algorithm, we talk about the localized dimension reduction procedures that are possible with these simple splits, and describe the construction of splitting dimensions $\Gras{D}$.

\subsection{Optimal splitting dimensions \texorpdfstring{$\Gras{D}^{*}$}{D*}}

Suppose that we found $\Gras{x}_{\ell}^{*}$ to split a leaf $\ell$ in all the dimensions $\{1,\ldots,d\}$. Before effectively splitting the leaf, we will check if, by chance, some dimensions can be removed from the splitting dimensions. Indeed, we will choose the splitting dimensions $\Gras{D}$ based on a statistical test whose hypothesis writes:

\begin{hypothesis}[Sparsity hypothesis $\mathcal{H}_j$]\label{hyp:dim_red1} Denoting $\Gras{U} \sim C$,  define for a given dimension $j \in \{1,\ldots,d\}$ the sparsity hypothesis as:
	$$\mathcal{H}_j: \quad\left(U_j  \indep U_{-j}\right) \,\vert \;\Gras{U} \in \ell\text{ and }U_j \,\vert \; \Gras{U} \in \ell \sim \mathcal{U}(\ell_j)\text{,}$$
	
	where we denoted $U_{j}$ the $j^{\text{th}}$ marginal of the random vector $U$ and $U_{-j}$ the vector of all other marginals. 
\end{hypothesis}

The hypothesis $\mathcal{H}_j$ literally says that, inside the leaf $\ell$, the dimension $j$ of the data is uniform and independent of the others. When $\mathcal{H}_j$ is accepted in a leaf $\ell$, we will reduce the dimension of the breakpoint $\Gras{x}_{\ell}^{*}$ in this leaf and in all child leaves accordingly, by removing $j$ from the set of splitting dimensions $\Gras{D}^{*}$. This will have the result that, inside the leaf $\ell$, the final model will consider the dimension $j$ to be uniform and independent of others.

To check this hypothesis using the integrated square error, as analyzed by~\cite{bowman1992}, suppose without loss of generality that $\ell$ is rescaled to $\mathbb{I}$, containing $n$ observations of the random vector $\Gras{U} \sim F$, for $F$ the restriction of $C$ to $\ell$, rescaled to $\mathbb{I}$ (note that $F$ is not a copula). This removes the conditioning in the hypothesis. Then, define the test statistic as follows.

\begin{definition}[Test statistic]
	Denote by $f_{\Gras{\mathrm{f}},\L}^{(n)}$ the piecewise constant density that will be estimated on data $\Gras{U} \sim F$, using the surrogate loss \eqref{eq:surrogate_loss_last_version}, and by $e_{j,n}(x) = \Esper{f_{\Gras{\mathrm{f}},\L}^{(n)}(x) | \mathcal{H}_j}$  the expectation of this estimate under $\mathcal{H}_j$. The test statistic is given by: 
	
	$$\mathcal{I}_j = \norme{e_{j,n} - f_{\Gras{\mathrm{f}},\L}^{(n)}},$$ 
	
	where $\L$, $e_{j,n}$ and $f_{\Gras{\mathrm{f}},\L}^{(n)}$ are stochastic objects, depending on the independent and identically distributed ($i.i.d.$) random vectors $\Gras{U}_1,\ldots,\Gras{U}_n$.	 
\end{definition}

This test statistic does not test hypothesis $\mathcal{H}_j$ \Emph{per se}, but rather tests $\mathcal{H}_j$ under the hypothesis of piecewise constant density. This is a weaker assertion, but it is enough to decide if local dimension reduction is possible or not given the current estimation stage and the data. More classical tests for independence and/or uniformity can be founded in~\cite{genest2019,yang2017,yao2018,dewet1980,song2009,hofert2014}.

The statistic $\mathcal{I}_j$ has the nice property that it is constructed as a square distance, and hence is always non-negative and is $0$ only under $\mathcal{H}_j$. On the other hand, it requires that we compute the full patchwork distribution in the two cases (under $\mathcal{H}_j$ or not), which can be costly. Instead, we can only compute \Emph{the next split}, which will reduce the computation and simplify the statistic. The drawback is that the test is weakened. The next property gives an empirical form of this statistic, using this simplification.

\begin{property}[Empirical form of the statistic]
	On a sample of data $(u_{i,j})_{n\times d}$, we can approximate the statistic $\mathcal{I}_j$ by: 
	
	\begin{equation}
	\label{eq:dim_red_stat}
	\widehat{\mathcal{I}_j} = \sum\limits_{k \in \L_{x^{*},\{1,\ldots,d\}\setminus \{j\}}} \left\{ \frac{\mathrm{f}_k^2}{\lambda(k)} + \sum\limits_{\substack{\ell \in \L_{x^{*},\{1,..,d\}}\\ \ell \subset k}} \left(\frac{\mathrm{f}_\ell^2}{\lambda(\ell)} - 2 \frac{\mathrm{f}_k\mathrm{f}_\ell}{\lambda(k)}\right)\right\}.
	\end{equation}
	
	\begin{proof} The estimator will cut on the same breakpoints on dimensions other than $j$ whether we work under $\mathcal{H}_j$ or not. This gives the definition of the partitions, and hence the expression for $e_{j,n}$ and $f_{\Gras{\mathrm{f}},\L}^{(n)}$, giving the expression of $\widehat{\mathcal{I}_j}$.
	\end{proof}
\end{property}

The law of the statistic \eqref{eq:dim_red_stat} under $\mathcal{H}_j$ cannot be computed explicitly. We use a Monte-Carlo simulation to compute the p-value of the test. To that purpose, simulate $n$ uniform random variables to replace $\Gras{u}_{.,j}$, and recompute the statistic \eqref{eq:dim_red_stat}, $T$ times. Indeed, note that under the null, the values of $\Gras{u}_{.,-j}$ can be held fixed.

The full localized dimension reduction procedure is formalized in Algorithm~\ref{algo:dim_red}.

\begin{algorithm}[H]
	\SetAlgoLined
	\KwData{$T \in \mathbb{N}$, a leaf $\ell$, observations $\Gras{u}_1,..,\Gras{u}_n \in \mathbb{\ell}$, and a threshold probability $\alpha$}
	\KwResult{The splitting dimensions $\Gras{D}^{*}$}
	Obtain $\Gras{x}_ {\ell}^{*}$ and $\L_{\Gras{x}_{\ell}^{*},\{1,..,d\}}$ by greedily minimizing the loss in \eqref{eq:surrogate_loss_last_version} on $\ell$.\\
	\ForEach{$j \in 1,\ldots,d$}{
		Denote by $s$ the statistic $\widehat{\mathcal{I}_j}$ given by \eqref{eq:dim_red_stat}.\\
		\ForEach{$i \in 1,\ldots,T$}{
			Simulate $n$ uniform random variables on $\ell_j$ and replace the $j$th coordinate of the data by this simulation\;
			Denote $s_i$ the value of the statistic from \eqref{eq:dim_red_stat} on this new dataset\;
		}
		Set $p_j = \frac{1}{T} \sum\limits_{i=1}^{T} \mathbb{1}_{s < s_i}$\\
	}
	\Return{$\Gras{D}^{*} = \left\{j \in 1,\ldots,d: p_j > 1 - \alpha\right\}$}\\
	\caption{Localized dimension reduction.\label{algo:dim_red}}
\end{algorithm} 

Note that the threshold $\alpha$ is a probability threshold for each individual test, and not a global threshold. To avoid multiple testing issues, it should be treated as a hyper-parameter and not a type $1$ error. We now formalize in the next subsection the complete estimation procedure.

\subsection{Full estimation procedure}

Suppose that we start from the independence copula, which writes $C_{\{1\},\{\mathbb{I}\}}$ in our framework. Then, if the sample of observations belonging to the sole leaf is too far from independence, i.e., if $\mathcal{H}_j$ does not hold for all $j$, we construct the first breakpoint by greedily minimizing the loss \eqref{eq:surrogate_loss_last_version} over the splitting point $\Gras{x}$. Rescaling the new leaves to $\mathbb{I}$ allows starting over and split again, until a proper stopping condition is reached: either there are no points anymore in the leaf or the leaf passes the uniformity tests. A third stopping condition is that the loss is no more reduced by splitting. 

To fasten the computation, we decided to ignore the copula constraints while splitting, and enforce them at the end on the constructed partition to correct the empirical weights. Later properties of convergence back up this decision, and this furthermore allows to parallelize the splitting process. Experiments showed that the algorithm that enforces the constraints at each split is much slower (since for the optimization of the breakpoint, sub-optimization corresponding to the quadratic program of Proposition~\ref{prop:final_quadratic_program} must be run for each evaluation) and does not provide much better results. 

More formally, Algorithm~\ref{algo:CORT} below states the complete estimation procedure.

\begin{algorithm}[H]
	\SetAlgoLined
	\KwData{Observed ranks $\Gras{u}_1,..,\Gras{u}_n \in \mathbb{I}$}
	\KwResult{Parameters $\Gras{p}$ and $\L$ of the estimated piecewise linear copula}
	Initialize the tree by $\L = \{\mathbb{I}\}$ and $\Gras{p}_\L = \{1\}$.\\
	\While{there exist leaves $\ell \in \L$ that are still splitable}{
		\ForEach{$\ell \in \L$ that is splitable}{
			Run Algorithm~\ref{algo:dim_red} in $\ell$ to find $\Gras{D}^{*}$.\\
			\If{$\Gras{D}^{*} \neq \emptyset$}{
				Find $\Gras{x}_{\ell,\Gras{D}^{*}}^{*}$ minimizing the surrogate loss \eqref{eq:surrogate_loss_last_version}.\\
				Set $\L = \L_{\Gras{x}_{\ell,\Gras{D}^{*}}^{*},\Gras{D}^{*}}$.
			}
		}
	}
	Compute $\Gras{p}_{\L}^{*} = P_{\L,\mathcal{C}_{L}}(\Gras{\mathrm{f}}_\L)$ from Proposition~\ref{prop:final_quadratic_program}.\\
	\Return{$(\Gras{p}_{\L}^{*},\L)$.}
	\caption{CORT estimation.\label{algo:CORT}}
\end{algorithm}

The resulting estimator of the copula density, denoted by $c_{\Gras{p}_{\L}^{*},\L}^{(n)}$, is called CORT for Copula Recursive Tree. The current implementation of this algorithm as an \textsf{R} package \texttt{cort} is available on \texttt{CRAN}. Note that the conditions for a leaf to be splitable can vary: by default, we consider that a leaf becomes non-splitable when it contains less than two observations. Options exist in the implementation to randomize the splitting dimensions (instead of optimizing them), or to constraint the maximum number of splitting dimensions to be inferior to some threshold. We do not provide the analysis of these parametrizations.

\begin{remark}
	Before step $11$ of Algorithm~\ref{algo:CORT}, the built model is simply the density estimate $f_{\Gras{\mathrm{f}}_{\L},\L}^{(n)}$. This additional outcome of the fitting procedure will be used in the next sections for performance analysis.
\end{remark}

\begin{remark}
	If we restrict the breakpoint possibility to all points with coordinates in $\left(\frac{i}{m+1}\right)_{i \in \{1,\ldots,m\}}$ where $m$ divides the number of observations $n$, this gives us only $m^d$ candidates. This corresponds to a form of checkerboard copula, see~\cite{cuberos2019} for more details. Since the ISE loss we use is tractable enough, the breakpoint can be chosen by directly minimizing the criterion over the continuous space $\mathbb{I}$ if the dimension of the problem is not too big.
\end{remark}

After discussing the consistency of this estimation procedure, Section~\ref{sec:4} discusses potential extensions, mainly through bagging principles. 

\subsection{Consistency}

We will show the consistency of the CORT estimator in the $\mathrm{L}_{2}$ almost-sure sense. Recall from \cite[Theorem 1]{ram2011} the following result about the \Emph{unconstrained} estimator. 

\begin{property}[Consistency of $f_{\Gras{\mathrm{f}}_{\L},\L}^{(n)}$]\label{prop:assymptotic_consistency_for_the_unconstraint_estimator}
	
	Assuming the maximum diameter of leaves goes to $0$ as $n$ goes to $\infty$, we have: $$\Proba{\lim\limits_{n \mapsto +\infty} \norme{f_{\Gras{\mathrm{f}}_{\L},\L}^{(n)} - c}=0}=1$$
	
\end{property}

A detailed proof, based on a generalization by Vapnik-Chervonenkis~\cite{vapnik2015} of the Glivenko-Cantelli Theorem, can be found in~\cite{ram2011}. Denote now by $\Gras{q}_{\L}$ the volumes given by the true copula on the leaves: $$\forall \ell \in \L ,\;q_\ell = \int\limits_{\ell} c(\Gras{u}) d\Gras{u}.$$ Then, one can easily check that $\Gras{q}_{\L} \in \mathcal{C}_{\L}$, and Proposition~\ref{prop:assymptotic_consistency_for_the_unconstraint_estimator} leads to the following useful corollary: 

\begin{corollary}\label{corr:dist_to_zero}
	$\distL{\Gras{\mathrm{f}}_{\L}}{\Gras{q}_{\L}} \rightarrow 0\text{, a.s}.$
\end{corollary}

Indeed, replacing $c$ by a piecewise constant density with partition $\L$ and weights $\Gras{q}_{\L}$ does not change the value of $\norme{f_{\Gras{\mathrm{f}}_{\L},\L}^{(n)} - c}$ and hence $\norme{f_{\Gras{\mathrm{f}}_{\L},\L}^{(n)} - c} = \distL{\Gras{\mathrm{f}}_{\L}}{\Gras{q}_{\L}}$.

\begin{definition}[Integrated constraint influence] Define the integrated constraint influence as the following squared norm: 
	
	\begin{equation}\label{eq:integrated_constraint_influence}
	\norme{c_{\Gras{p}_{\L}^{*},\L}^{(n)} - f_{\Gras{\mathrm{f}}_{\L},\L}^{(n)}} = \distL{\Gras{p}_{\L}^{*}}{\Gras{\mathrm{f}}_{\L}}
	\end{equation}
	
\end{definition}

In the simulation studies in Section~\ref{sec:5}, this quantity will be monitored via burn-in techniques, to see how it behaves as $n$ grows. Proposition~\ref{prop:asymptotic_constraint_influence_is_zero} below gives the corresponding theoretical result.

\begin{property}[Asymptotic effect of constraints]\label{prop:asymptotic_constraint_influence_is_zero} As $n \rightarrow \infty$, the integrated constraint influence goes to $0$. 
	
	\begin{proof} Recall from Proposition~\ref{prop:final_quadratic_program} that $\Gras{p}_{\L}^{}*$ is the orthogonal projection of $\Gras{\mathrm{f}}_{\L}$ into $\mathcal{C}_{\L}$. Since $\Gras{q}_{\L} \in \mathcal{C}_{\L}$, we have that $\distL{\Gras{\mathrm{f}}_{\L}}{\Gras{p}_{\L}^{*}}\le\distL{\Gras{f}_{\L}}{\Gras{q}_{\L}}$, which concludes the argument by Corollary~\ref{corr:dist_to_zero}.
	\end{proof}
\end{property}	The consistency of the estimator is now easy to obtain.

\begin{property}[Consistency] For $c$ the density of the true copula, assuming the diameter of the leaves goes to $0$ as $n$ goes to $\infty$, the estimator $c_{\Gras{p}_{\L}^{*},\L}^{(n)}$ is consistent, i.e.: 
	$$\Proba{\lim\limits_{n \mapsto +\infty} \norme{c_{\Gras{p}_{\L}^{*},\L}^{(n)} - c}=0}=1$$	
	
	\begin{proof}	Remark that $\norme{c_{\Gras{p}_{\L}^{*},\L}^{(n)} - c} = \distL{\Gras{p}_{\L}^{*}}{\Gras{q}_{\L}}$. Then, $\Gras{p}_{\L}^{*}$ being the orthogonal projection of $\Gras{\mathrm{f}}_{\L}$ into $\mathcal{C}_{\L}$, and $\Gras{q}_{\L}$ being in $\mathcal{C}_{\L}$, we have $\distL{\Gras{p}_{\L}^{*}}{\Gras{q}_{\L}} \le \distL{\Gras{\mathrm{f}}_{\L}}{\Gras{q}_{\L}}$ and we conclude by Corollary~\ref{corr:dist_to_zero}.			
	\end{proof}
\end{property}

%In the next section, we consider bagging methods to increase the performance of the model. 

\section{Bagging of density estimators\label{sec:4}}

While using the CORT algorithm, overfitting is probable since we do not use any kind of pruning. Although an implementation of pruning might be possible in our setting (remove one leaf, compute new weights through the adequate quadratic program, and compare the loss in the objective function to the final loss, then design a criterion to know where to stop), we preferred to augment the 'complexity' of the model by mixing (bagging) rather than reduce it by pruning, since the model is already quite simple.
Furthermore, with a minimum node size (minimum number of points in a leaf to initiate a split) set to one, the model produces an average number of leaves growing exponentially with the dimension, a lot of them being empty. The bagging procedure we propose in this section should overcome these two issues by mixing different trees fitted on different resamples of the dataset. 

Two estimates $\hat{g}_1$ and $\hat{g}_2$ of a function $g$ can be bagged into an estimate $\frac{\hat{g}_1 + \hat{g}_2}{2}$, candidate for the estimation of $g$ if the estimates are close to be uncorrelated. This principle gave rise to the bagging algorithm in regression, developed by Breiman in~\cite{breiman1996}. In density estimation, bagging can also be exploited: if the estimator has a high variance and a small bias, then bagging it might yield a better result, regarding the bias/variance trade-off.

Recall that, while bootstrapping over $n$ observations, the chance that an observation does not appear in the bootstrap sample is given by $\left(1- \frac{1}{n}\right)^n \rightarrow \frac{1}{e}$. Hence, asymptotically $36.8$ percent of the dataset will end up \Emph{out-of-bag}. These samples can be used to check for accuracy of the model, and even set hyper-parameters when some are needed. Following the work of~\cite{breiman1996} in regression models,~\cite{sain1994} formalized the cross-validation and bagging process for density estimation. Note that usually, the \Emph{leave-one-out} method is used in kernel density estimation to select hyper-parameters (mainly the bandwidth), see, e.g.,~\cite{miecznikowski2010}, but the more involved \Emph{out-of-bag} procedure we propose is inspired by~\cite{wu2018}.

In the following, we denote by $\hat{C}$ the empirical copula of the whole dataset $\Gras{u}$, a $n\times d$ matrix.

\begin{definition}[CORT Forest]\label{def:CORT_forest}
	Define $\Gras{u}^{(1)},\ldots,\Gras{u}^{(N)}$ as $N$ bootstrap resamples of the same size $n$, and $\hat{c}^{(1)}, ...,\hat{c}^{(N)}$ (resp. $\hat{C}^{(1)}, ...,\hat{C}^{(N)}$) the densities (resp. d.f.) of the CORT estimators on these resamples fitted by Algorithm~\ref{algo:CORT}. Define the CORT forest with weights $\Gras{\omega} = (\omega_1,\ldots,\omega_N)$ as the mixture distribution with density: $$\hat{c}_{\Gras{\omega}}(\Gras{v}) = \sum_{j=1}^{N} \omega_j  \hat{c}^{(j)}(\Gras{v}) \; \forall \Gras{v} \in \mathbb{I}.$$
	
\end{definition}

For each observation $\Gras{u}_i$ in the original training set, we recall the out-of-bag density estimate as: $$\hat{c}_{\Gras{\omega}}^{oob}(\Gras{u}_i) = \frac{\sum\limits_{j=1}^{N}  \omega_j \hat{c}^{(j)}(\Gras{u}_i)\mathbb{1}_{\Gras{u}_i \notin \Gras{u}^{(j)}}}{\sum\limits_{j=1}^{N} \omega_j \mathbb{1}_{\Gras{u}_i \notin \Gras{u}^{(j)}}}.$$

We define $\hat{C}_{\Gras{\omega}}$ and $\hat{C}_{\Gras{\omega}}^{oob}$ accordingly. Note that $\hat{c}_{\Gras{\omega}}^{oob}$ is not a proper density as it may fail to sum to $1$ and it is only defined on the observation points. Because trees were fitted independently, this is however, on observation points, a good approximation of the forest density itself. From $\hat{c}_{\Gras{\omega}}^{oob}$, based on~\cite{miecznikowski2010},~\cite{wu2018} defines out-of-bag version of common fitting statistics as follows.

\begin{definition}[Out-of-bag statistics]\label{def:oob_stats} Define respectively the out-of-bag empirical integrated square error, the out-of-bag Kullback-Leibler divergence and two out-of-bag Cramer-Von-Mises distances associated to the forest as: 
	\begin{align*}
	\hat{J}(\hat{c}_{\Gras{\omega}}) &=  \norme{\hat{c}_{\Gras{\omega}}} - \frac{2}{n} \sum\limits_{i=1}^n \hat{c}_{\Gras{\omega}}^{oob}(\Gras{u}_i),\qquad\hat{K}(\hat{c}_{\Gras{\omega}}) =  -\frac{1}{n}\sum\limits_{i=1}^n\ln\left(\hat{c}_{\Gras{\omega}}^{oob}(\Gras{u}_i)\right)\\
	\hat{M}(\hat{C}_{\Gras{\omega}}) & = \frac{1}{n} \sum_{i=1}^{n} \left(\hat{C}_{\Gras{\omega}}^{oob}(\Gras{u}_i) - \hat{C}(\Gras{u}_i)\right)^2,\qquad\hat{N}(\hat{C}_{\Gras{\omega}}) = \frac{1}{n} \sum_{i=1}^{n} \hat{C}_{\Gras{\omega}}(u_i)^2 - 2 \hat{C}_{\Gras{\omega}}^{oob}(\Gras{u}_i) \hat{C}(\Gras{u}_i).
	\end{align*}
\end{definition}

Note that $\hat{K}(\hat{c}_{\Gras{\omega}})$ is obtained as an empirical version of $\int \ln\left(\frac{c(\Gras{u})}{\hat{c}_{\Gras{\omega}}^{oob}(\Gras{u})}\right)dC(\Gras{u})$. $\hat{M}(\hat{c}_{\Gras{\omega}})$ estimates the Cramer-Von-Mises distance between the out-of-bag d.f. of the forest and the empirical copula. On the other hand, $\hat{N}(\hat{c}_{\Gras{\omega}})$ keeps the true norm of the model, in the same spirit as $\hat{J}(\hat{c}_{\Gras{\omega}})$, see~\cite{miecznikowski2010,wu2018} for more details about these cross-validation tools.

We denote: $$\Omega = \left\{\omega \in \mathbb R^N,\; \sum\limits_{j=0}^N \omega_j = 1 \text{ and } \omega_j\ge 0\; \forall j \in \{1,\ldots,N\} \right\}$$ the set of possible weights. Indeed, we want $\hat{c}_{\Gras{\omega}}$ to be a density, which implies that weights sum to one and are positives. We call \Emph{optimal} the forest with the weights minimizing $\hat{J}(\hat{c}_{\Gras{\omega}})$. We use an optimization program to find weights given the resampling and the trees: $$\Gras{\omega}^{*} = \argmin\limits_{\Gras{\omega} \in \Omega}\;\hat{J}(\hat{c}_{\Gras{\omega}}).$$ Note that this method makes out-of-bag observations contribute to the final estimation through weights. 

The forest estimation is studied in the next section. Note that the constructions of this section are the same if you use another base estimator than the CORT estimator: in the next section, we will use these tools to compare bagging of the CORT estimator with bagging of other copula estimators.

\section{Investigation of performance\label{sec:5}}

In this section, we investigate the performance of the proposed estimation procedure on several simulated datasets. To ensure reproducibility, we provide all the datasets in the \textsf{R} package \texttt{cort}, with the code and parameters needed to re-simulate them. We will compare our results with several other models: \begin{itemize}
	\item Deheuvel's empirical copula, hereafter denoted by ``Empirical'';
	\item The empirical beta copula~\cite{segers2017}, hereafter denoted ``Beta'';
	\item A Checkerboard copula~\cite{cuberos2019} with $m=10$, denoted by ``Cb(m=10)'';
	\item Another less-precise Checkerboard copula with $m=5$, denoted by ``Cb(m=5)''.
\end{itemize}

Recall from~\cite{cuberos2019} that a checkerboard copula with a parameter $m$ is a piecewise linear copula with a partition $\L$ composed of a regular grid of hypercubes with side length $m^{-1}$. See the reference for details about the empirical beta copula. They are all non-parametric or semi-parametric models, with a straightforward estimation procedure. The two checkerboard copulas are provided by the \textsf{R} package \texttt{cort}, and the empirical beta copula, as well as the empirical copula, are from the \textsf{R} package \texttt{copula}. 

We will compare results of different models in terms of dependence measures, Kendall tau and Spearman rho first. Then, we will look at predictive performance metrics defined in the previous section, $\hat{J}$, $\hat{K}, \hat{M}$ and $\hat{N}$, computed via a weighted bagging of each model. 

\begin{dataset}[Piecewise linear copula]\label{data:recoveryourself}
	This dataset is a simple test: we simulate random samples from a density \Emph{inside} the piecewise copula class, and test whether the estimator can recover it. For that, we will use a $2$-dimensional sample with $500$ observations, uniform on $\mathbb{I}$, and apply the following function: $$h_1(\Gras{u}) = \left(u_1, \frac{u_2 + \mathbb{1}_{u_1 \le \frac{1}{4}} + 2\mathbb{1}_{u_1 \le \frac{1}{2}} + \mathbb{1}_{\frac{3}{4} \le u_1}}{4}\right).$$
\end{dataset}
% \begin{figure}[H]
% 	\centering
% 	\includegraphics[width=0.40\textwidth]{recoveryourself_data.png}
% 	\caption{Pairs-plot of Dataset~\ref{data:recoveryourself}: this bivariate dataset has already a piecewise linear structure.}
% \end{figure}
% % % % % % % % % % % % % % % % % % % % % % % % 

\begin{figure}[H]
	\centering
	\includegraphics[width=0.7\textwidth]{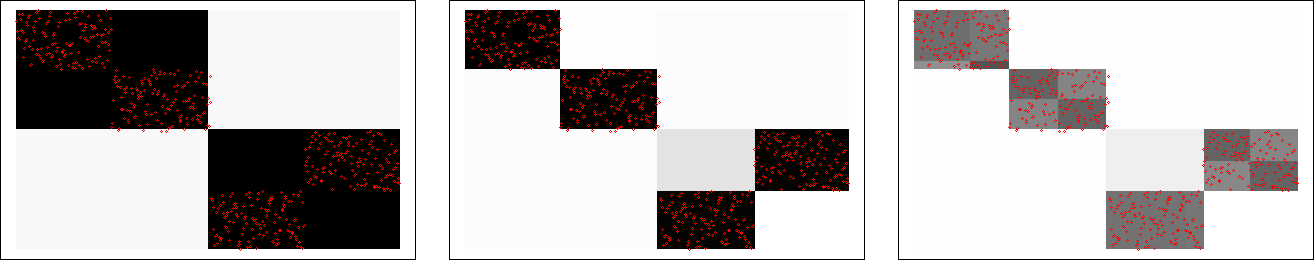}
	\caption{(Dataset~\ref{data:recoveryourself}) Running example. Data points are in red, and darker zones mean boxes with higher weights (black is the maximum possible in each leaf). On the left, the model after the first split. In the middle, the model after the next round of splits, and on the right after the third round.}
	\label{fig:recoveryourself_running_exemple}
\end{figure}

\begin{figure}[H]
	\centering
	\subfloat[]{{\includegraphics[width=0.40\textwidth]{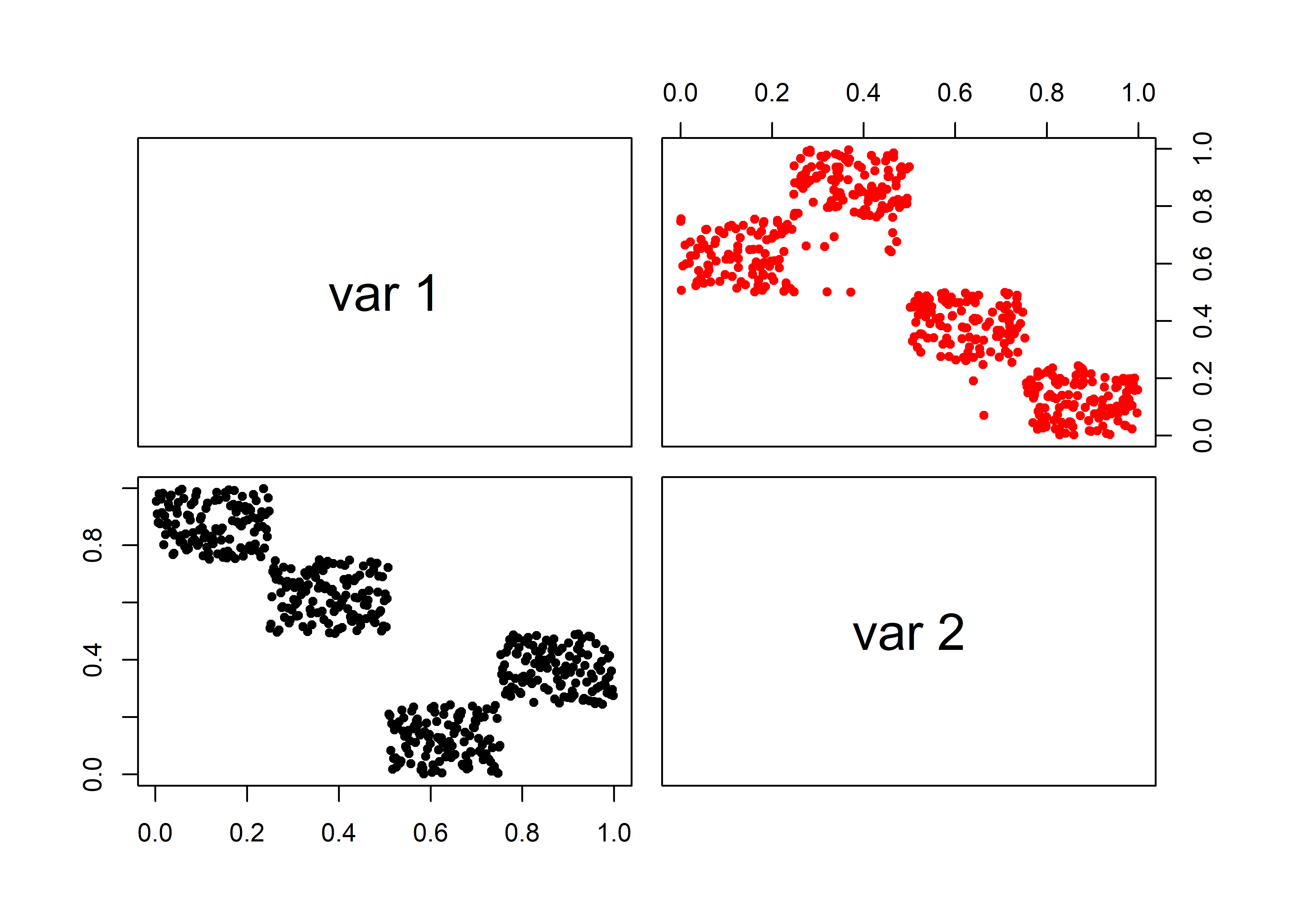}}}%
	\qquad
	\subfloat[]{{\includegraphics[width=0.40\textwidth]{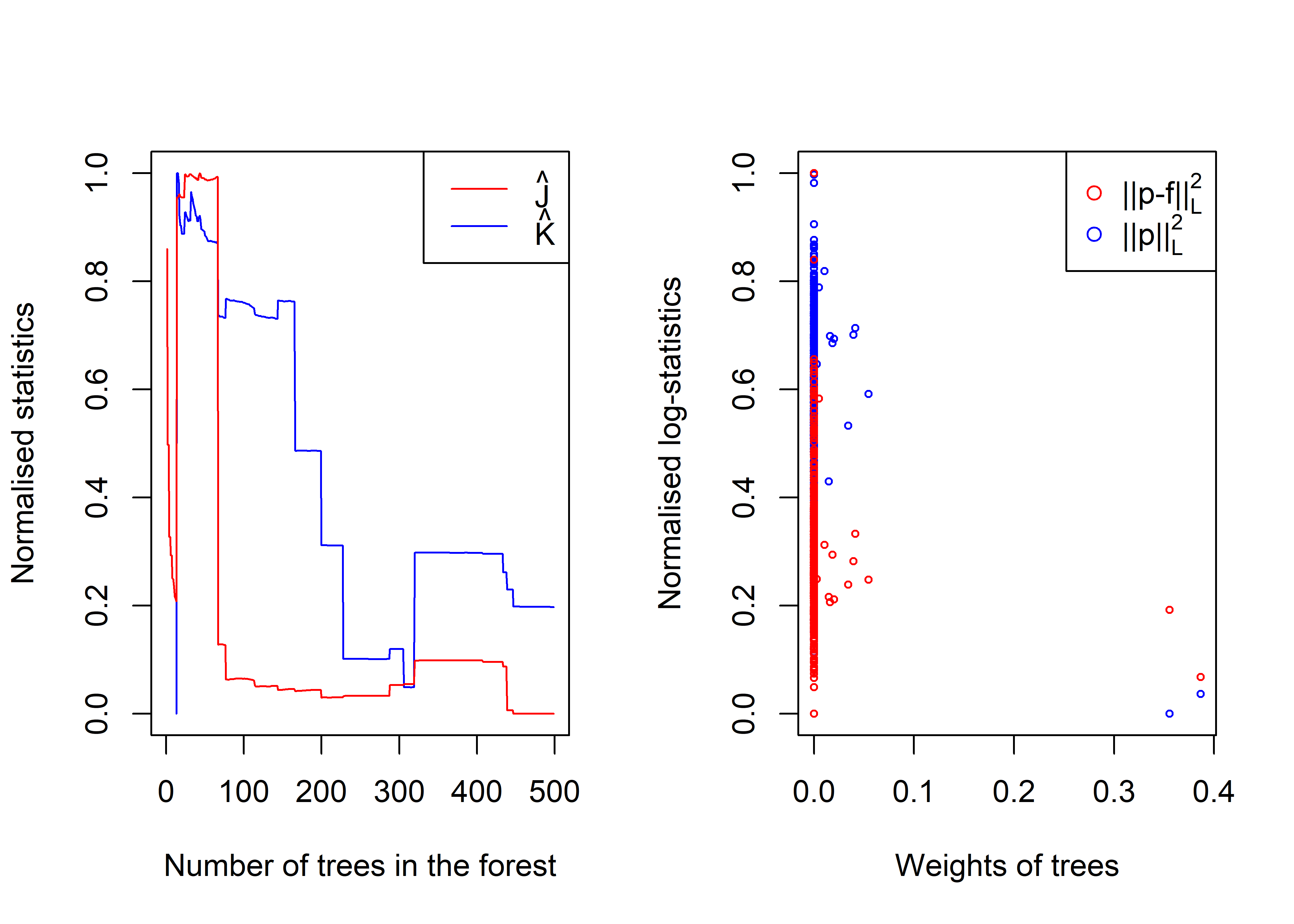} }}%
	\caption{(Dataset~\ref{data:recoveryourself})  (a) The CORT estimator: in black, lower left, the input data. In red, upper-right, a simulation from the estimated tree. (b) Statistics of the forest: on the left, $\hat{K}$ and $\hat{J}$ in function of the number of trees. On the right, the Integrated Constraint Influence and square norm of each tree against the weight of the tree in the forest.}%
	\label{fig:recoveryourself}%
\end{figure}

To illustrate the behavior of the algorithm, we propose in Fig.~\ref{fig:recoveryourself_running_exemple} to look at a running example of the first few splits. Fig.~\ref{fig:recoveryourself} (a) shows a simulation from the final (fully fitted) CORT estimator on Dataset~\ref{data:recoveryourself}, along with results from the CORT forest.

%We see that the model identified correctly the first splitting points. 

% \begin{figure}[H]
% 	\centering
% 	\subfloat[In black, lower left, the input data. In red, upper-right, a simulation from the estimated tree.]{{\includegraphics[width=0.40\textwidth]{recoveryourself_tree.png}}}%
% 	\qquad
% 	\subfloat[On the left, $\hat{K}$ and $\hat{J}$ in function of the number of trees. On the right, the Integrated Constraint Influence and square norm of each tree against the weight of the tree in the forest.]{{\includegraphics[width=0.40\textwidth]{recoveryourself_forest_stats.png} }}%
% 	\caption{(Dataset~\ref{data:recoveryourself})  (a) The CORT estimator (b) Statistics of the forest.}%
% 	\label{fig:recoveryourself}%
% \end{figure}

We observe on Fig.~\ref{fig:recoveryourself} (a) that the algorithm splits the space as requested. The few simulated points \Emph{outside} the four main boxes are there because the algorithm did not split \Emph{exactly} on $(\frac{1}{4},\frac{3}{4})$, $(\frac{1}{2},\frac{1}{2})$ and $(\frac{3}{4},\frac{1}{4})$, and the constraints forced him to put some weight on some leaves that do not contain points. 

Bagging the CORT algorithm on this dataset, we obtain statistics given by Fig.~\ref{fig:recoveryourself} (b). We observe that the out-of-bag statistics are decreasing in the number of trees fitted, although the weighting of the forest did select less than $10$ trees over $500$. Altogether, the algorithm succeeded into finding the right breakpoints. A comparison of the fit in terms of dependence measure to other models is available in Table~\ref{tab:recoveryourself_table_dep_mes}.

% \begin{table}[H]
% 	\caption{\label{tab:recoveryourself_table_dep_mes}Obtained dependence measures of several models on Dataset \ref{data:recoveryourself}. The first column is the goal, others are concurrent models.}
% 	\centering
% 	\begin{tabu} to \linewidth {>{\raggedright}X>{\raggedright}X>{\raggedright}X>{\raggedright}X>{\raggedright}X>{\raggedright}X>{\raggedright}X}
% 	\toprule
% 	  & Empirical & Cb(m=10) & Cb(m=5) & Beta & CORT & Bagged CORT\\
% 	\addlinespace[0.3em]
% 	\multicolumn{7}{l}{\Bold{Kendall Tau}}\\
% 	\hspace{1em}$\tau_{1,2}$ & -0.534 & -0.515 & -0.465 & -0.527 & -0.525 & -0.393\\
% 	\addlinespace[0.3em]
% 	\multicolumn{7}{l}{\Bold{Spearman Rho}}\\
% 	\hspace{1em}$\rho_{1,2}$ & -0.773 & -0.757 & -0.697 & -0.766 & -0.762 & -0.604\\
% 	\bottomrule
% 	\end{tabu}
% \end{table}

\begin{table}[H]
	\caption{\label{tab:recoveryourself_table_dep_mes}Obtained dependence measures (Kendall tau and Spearman rho) of several models on Dataset \ref{data:recoveryourself}. The first column is the goal, others are concurrent models.}
	\centering
	\begin{tabu} to \linewidth {>{\raggedright}X>{\raggedright}X>{\raggedright}X>{\raggedright}X>{\raggedright}X>{\raggedright}X>{\raggedright}X}
	\toprule
	  & Empirical & Cb(m=10) & Cb(m=5) & Beta & CORT & Bagged CORT\\
	\Bold{$\tau$} & -0.534 & -0.515 & -0.465 & -0.527 & -0.525 & -0.393\\
	\Bold{$\rho$} & -0.773 & -0.757 & -0.697 & -0.766 & -0.762 & -0.604\\
	\bottomrule
	\end{tabu}
\end{table}

On Table~\ref{tab:recoveryourself_table_dep_mes}, we display pairwise dependence measures (Kendall $\tau$ and Spearman $\rho$) of the obtained fits. To read these measures, consider the first column, corresponding to the empirical copula, as the goal for other models. We observe that all models perform correctly regarding dependence measures on this dataset, although the Checkerboard with $m=5$ (which has a pretty rough partition, not including the $\frac{1}{4}$ multiples) has a Spearman $\rho$ a little too high. Furthermore, bagging the CORT model gives the worst results.

Performing a standard weighted bagging procedure, we obtain fit statistics $\hat{K},\hat{J},\hat{M}$ and $\hat{N}$, displayed in Table~\ref{tab:recoveryourself_table_bagging_result}. The experiment fits every model $500$ times on resamples of the dataset and weights the resulting models to minimize $\hat{J}$. 

\begin{table}[H]

	\caption{\label{tab:recoveryourself_table_bagging_result}Results of the bagging of each model on Dataset \ref{data:recoveryourself}. Each row represents a different performance metric: in all cases, lower is better.}
	\centering
	\begin{tabu} to \linewidth {>{\raggedright}X>{\raggedright}X>{\raggedright}X>{\raggedright}X>{\raggedright}X>{\raggedright}X}
	\toprule
	  & Empirical & Cb(m=10) & Cb(m=5) & Beta & CORT\\
	$\hat{J}(\hat{c}_{\Gras{\omega}})$ & 0.002 & -3.16 & -2.37 & -2.98 & \Bold{-4.81}\\
	$\hat{K}(\hat{c}_{\Gras{\omega}})$ & Inf & -1.06 & -0.743 & \Bold{-1.18} & -0.837\\
	$\hat{M}(\hat{c}_{\Gras{\omega}})$ & \Bold{8.72e-06} & 4.13e-05 & 0.000281 & 1.58e-05 & 0.000633\\
	$\hat{N}(\hat{c}_{\Gras{\omega}})$ & \Bold{-0.0277} & -0.0277 & -0.0275 & -0.0277 & -0.0262\\
	\bottomrule
	\end{tabu}
\end{table}

Note that the Kullback-Leibler out-of-bag divergence $\hat{K}$ is infinite for the empirical copula, since it does not assign weights to points it did not see. We observe in bagging results that the predictive performance of the checkerboards is quite poor, that the forest based on empirical Beta copulas is a lot better and that the CORT forest generalized even better regarding the density-based measures $\hat{J}$ and $\hat{K}$. Note that the  bagging of empirical beta copula is a very powerful model.

\begin{dataset}[Another piecewise linear copula]\label{data:impossible}
	As for Dataset~\ref{data:recoveryourself}, we simulate from a density inside the piecewise linear copula class, by applying the function:$$h_2(\Gras{u}) = \left(u_1,          \frac{u_2}{2} + \frac{1}{2}\mathbb{1}_{u_1 \notin [\frac{1}{3}, \frac{2}{3})}        \right)$$ to a $200\times2$ uniform sample, and taking ranks.
\end{dataset}
% \begin{figure}[H]
% 	\centering
% 	\includegraphics[width=0.40\textwidth]{impossible_data.png}
% 	\caption{Pairs-plot of Dataset~\ref{data:impossible}: this dataset splits intentionally the space in a ternary way.}
% \end{figure}
% % % % % % % % % % % % % % % % % % % % % % % % 

\begin{figure}[H]
	\centering
	\subfloat[]{{\includegraphics[width=0.40\textwidth]{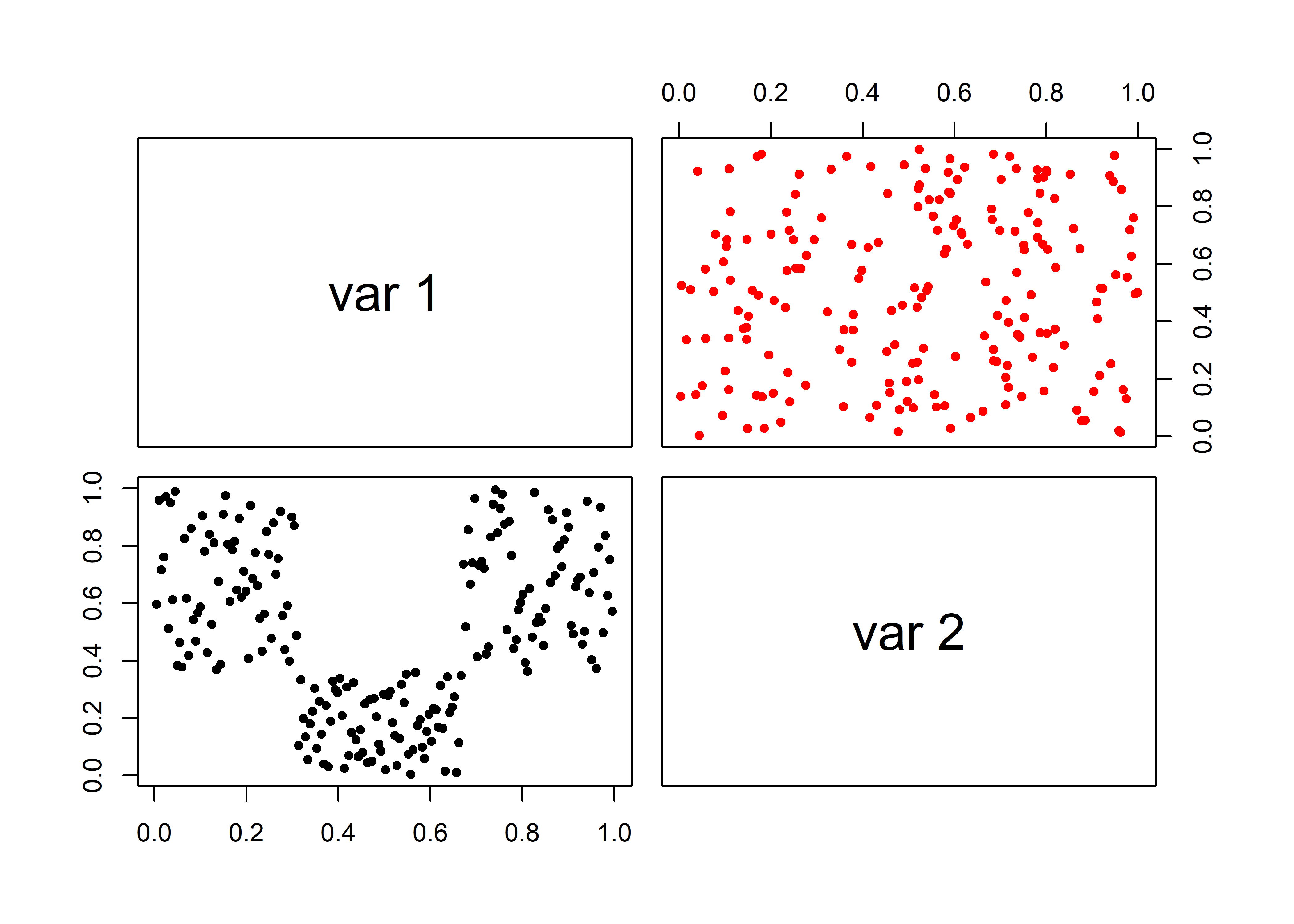}}}%
	\qquad
	\subfloat[]{{\includegraphics[width=0.40\textwidth]{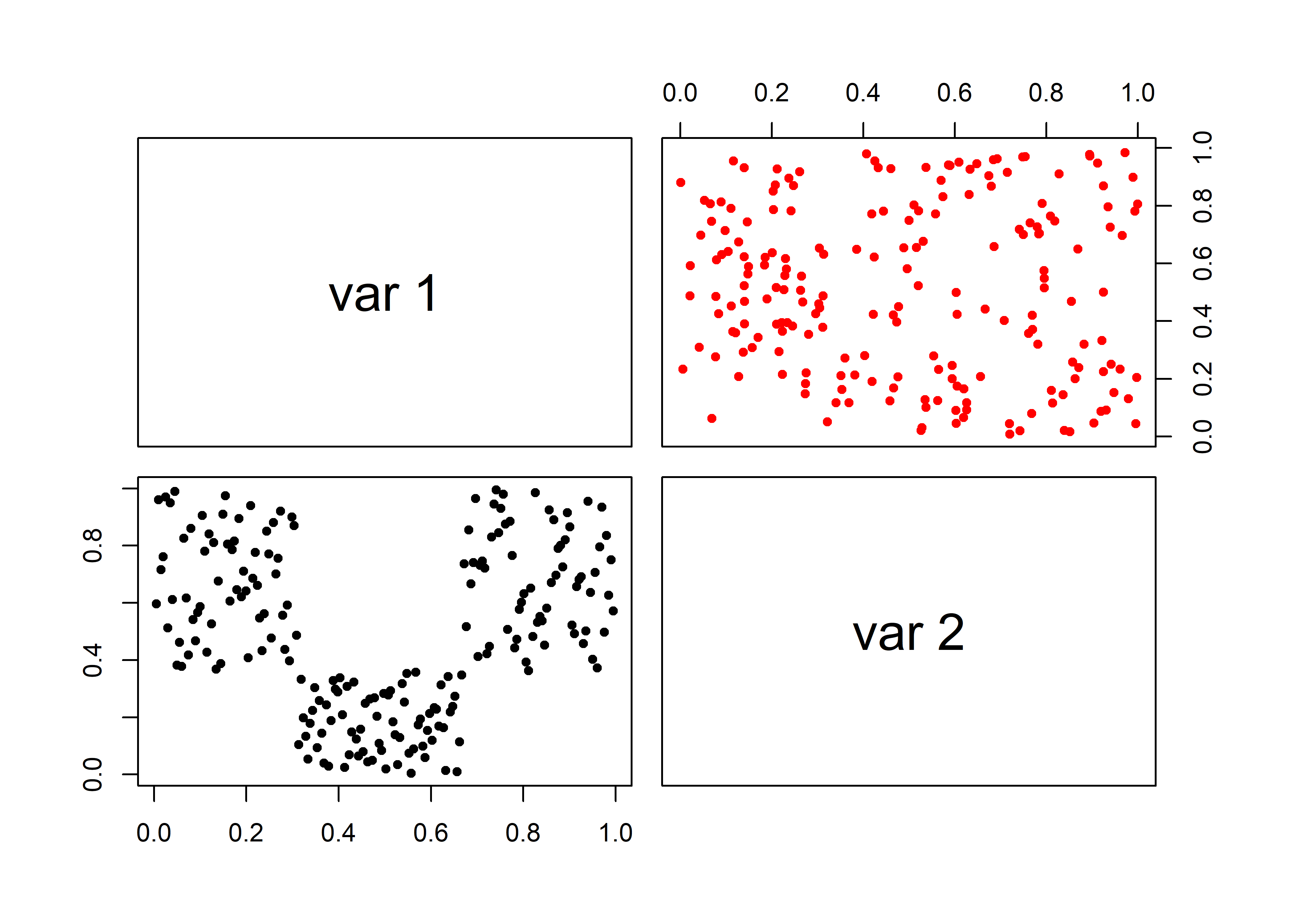} }}%
	\caption{(Dataset~\ref{data:impossible})  (a) The estimated tree: in black, lower left, the input data. In red, upper-right, a simulation from the estimated tree. (b) The estimated forest: in black, lower left, the input data. In red, upper-right, a simulation from the estimated forest.}%
	\label{fig:impossible}%
\end{figure}

This second dataset is also in the piecewise linear class, but it splits the space in a ternary way, which the recursive splitting procedure of the CORT estimator cannot reproduce. In Fig.~\ref{fig:impossible} (a), you can observe simulation from the CORT copula and the CORT forest fitted on Dataset~\ref{data:impossible}.

Remember that the algorithm splits recursively on only one breakpoint. If the data is split in a ternary way, as in Dataset~\ref{data:impossible}, it will not succeed. Looking at details from the fitting procedure, we see that the constraints forced the algorithm to return exactly the independence copula (by setting weights of each leaf equal to its volume). Indeed, while splitting, two splits in adjacent leaves are \Emph{not} synchronized: since the optimization routines are independent of each other, it is unlikely that they return the same breakpoint value for a given dimension, meaning that weights will not be transferable between the two zones: in our case, the constraints forced the weights back to independence. We see that the forest tried to correct this behavior, but the result is quite bad. 
%Table~\ref{tab:impossible_table_dep_mes} and Table~\ref{tab:impossible_table_bagging_result} contain the same results as for the previous dataset.
%
% \begin{table}[H]
% 	\caption{\label{tab:impossible_table_dep_mes}Obtained dependence measures of several models on Dataset \ref{data:impossible}.  The first column is the goal, others are concurrent models.}
% 	\centering
% 	\begin{tabu} to \linewidth {>{\raggedright}X>{\raggedright}X>{\raggedright}X>{\raggedright}X>{\raggedright}X>{\raggedright}X>{\raggedright}X}
% 	\toprule
% 	  & Empirical & Cb(m=10) & Cb(m=5) & Beta & CORT & Bagged CORT\\
% 	  \Bold{$\tau$} & 0.006 & 0.005 & 0.004 & 0.014 & 0 & -0.010\\
% 	  \Bold{$\rho$} & 0.010 & 0.011 & 0.009 & 0.021 & 0 & -0.014\\
% 	\bottomrule
% 	\end{tabu}
% \end{table}
%
The dependence measures ($\tau$ and $\rho$) are here structurally $0$, and every model respected this correctly. Table~\ref{tab:impossible_table_bagging_result}
% and Fig.~\ref{fig:impossible_Box plot} 
shows the same results as for the previous model.

\begin{table}[H]
	\caption{\label{tab:impossible_table_bagging_result}Results of the bagging of each model on Dataset \ref{data:impossible}. Each row represents a different performance metric: in all cases, lower is better.}
	\centering
	\begin{tabu} to \linewidth {>{\raggedright}X>{\raggedright}X>{\raggedright}X>{\raggedright}X>{\raggedright}X>{\raggedright}X}
	\toprule
	  & Empirical & Cb(m=10) & Cb(m=5) & Beta & CORT\\
	$\hat{J}(\hat{c}_{\Gras{\omega}})$ & 0.00501 & -1.48 & -1.39 & -1.21 & \Bold{-2.55}\\
	$\hat{K}(\hat{c}_{\Gras{\omega}})$ & Inf & -0.346 & -0.271 & \Bold{-0.426} & -0.311\\
	$\hat{M}(\hat{c}_{\Gras{\omega}})$ & \Bold{3.03e-05} & 5.16e-05 & 0.00016 & 7.38e-05 & 0.000907\\
	$\hat{N}(\hat{c}_{\Gras{\omega}})$ & \Bold{-0.134} & -0.134 & -0.134 & -0.134 & -0.134\\
	\bottomrule
	\end{tabu}
\end{table}

% \begin{figure}[H]
% 	\centering
% 	\subfloat[Box plot of log-normalized $\hat{P}$]{{\includegraphics[width=0.40\textwidth]{impossible_boxplotP_lognorm.png}}}%
% 	\qquad
% 	\subfloat[Box plot of log-normalized $\hat{Q}$]{{\includegraphics[width=0.40\textwidth]{impossible_boxplotQ_lognorm.png}}}%
% 	\caption{(Dataset~\ref{data:impossible})  Box plots of resulting errors for $20$ resamples for each model (the lower the better). $\hat{P}$ is focused on the distribution function and $\hat{Q}$ on the density.}%
% 	\label{fig:impossible_Box plot}%
% \end{figure}
% % and the box plot of $\hat{P}$ error confirms this analysis: in this case, the empirical beta copula performs a lot better.

Note that the out-of-bag ISE did not decrease during the fit of the forest, so that more complexity was added without any gain. We also observe that many trees in the forest have zero weights, they do not fit the data enough. Finally,  Table~\ref{tab:impossible_table_bagging_result} shows that the predictive performance of the CORT algorithm is quite bad.
% and the box plot of $\hat{P}$ error confirms this analysis: in this case, the empirical beta copula performs a lot better.

We now turn ourselves to the third dataset.

% % % % % % % % % % % % % % % % % % % % % % % % %  DATA CLAYTON
% \noindent\begin{minipage}{0.45\textwidth}
% 	\begin{dataset}[Modified Clayton]\label{data:clayton}
% 		This dataset is a simulation of $200$ points from a $3$-dimensional Clayton copula~\cite{joe2014} with $\theta = 7$ (hence highly dependent), for the first, third and fourth marginals. The second marginal is added as independent uniform draws. Lastly, the third marginal is flipped, inducing a negative dependence structure.
% 	\end{dataset}
% \end{minipage}
% \hfill%
% \begin{minipage}{0.45\textwidth}% adapt widths of minipages to your needs
% 	\centering\raisebox{\dimexpr \topskip-\height}{%
% 		\includegraphics[width=\textwidth]{clayton_data.png}}
% 	\captionof{figure}{Pairs-plot of Dataset~\ref{data:clayton}}
% \end{minipage}%% % % % % % % % % % % % % % % % % % % % % % % % 

% % % % % % % % % % % % % % % % % % % % % % % % DATA CLAYTON
\begin{dataset}[Modified Clayton]\label{data:clayton}
	This dataset is a simulation of $200$ points from a $3$-dimensional Clayton copula~\cite{joe2014} with $\theta = 7$ (hence highly dependent), for the first, third and fourth marginals. The second marginal is added as independent uniform draws. Lastly, the third marginal is flipped, inducing a negative dependence structure.
\end{dataset}
% \begin{figure}[H]
% 	\centering
% 	\includegraphics[width=0.40\textwidth]{clayton_data.png}
% 	\caption{Pairs-plot of Dataset~\ref{data:clayton}: this dataset is generated from a simple Clayton copula.}
% \end{figure}
% % % % % % % % % % % % % % % % % % % % % % % % 

\begin{figure}[H]
	\centering
	\subfloat[]{{\includegraphics[width=0.40\textwidth]{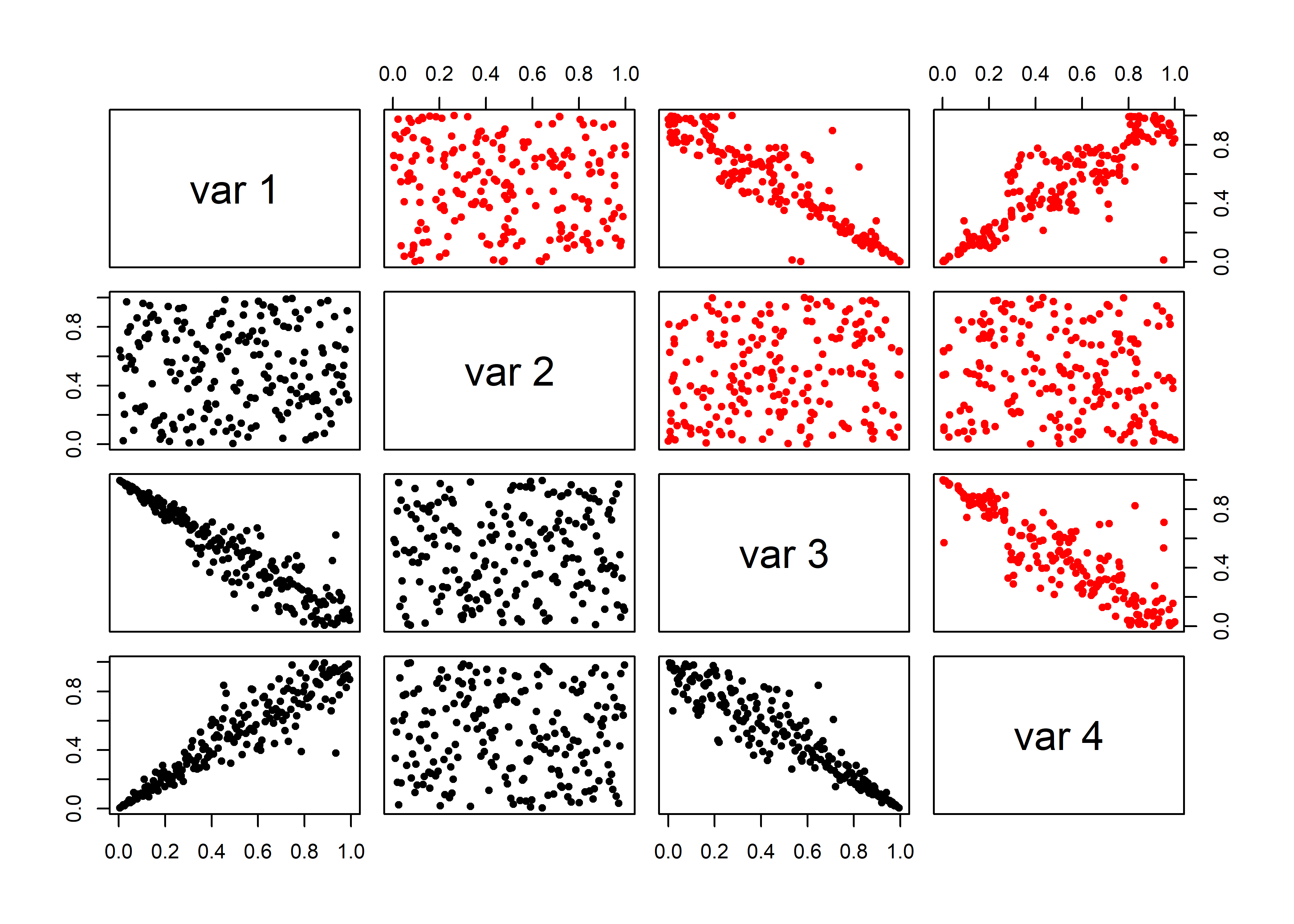} }}%
	\qquad
	\subfloat[]{{\includegraphics[width=0.50\textwidth]{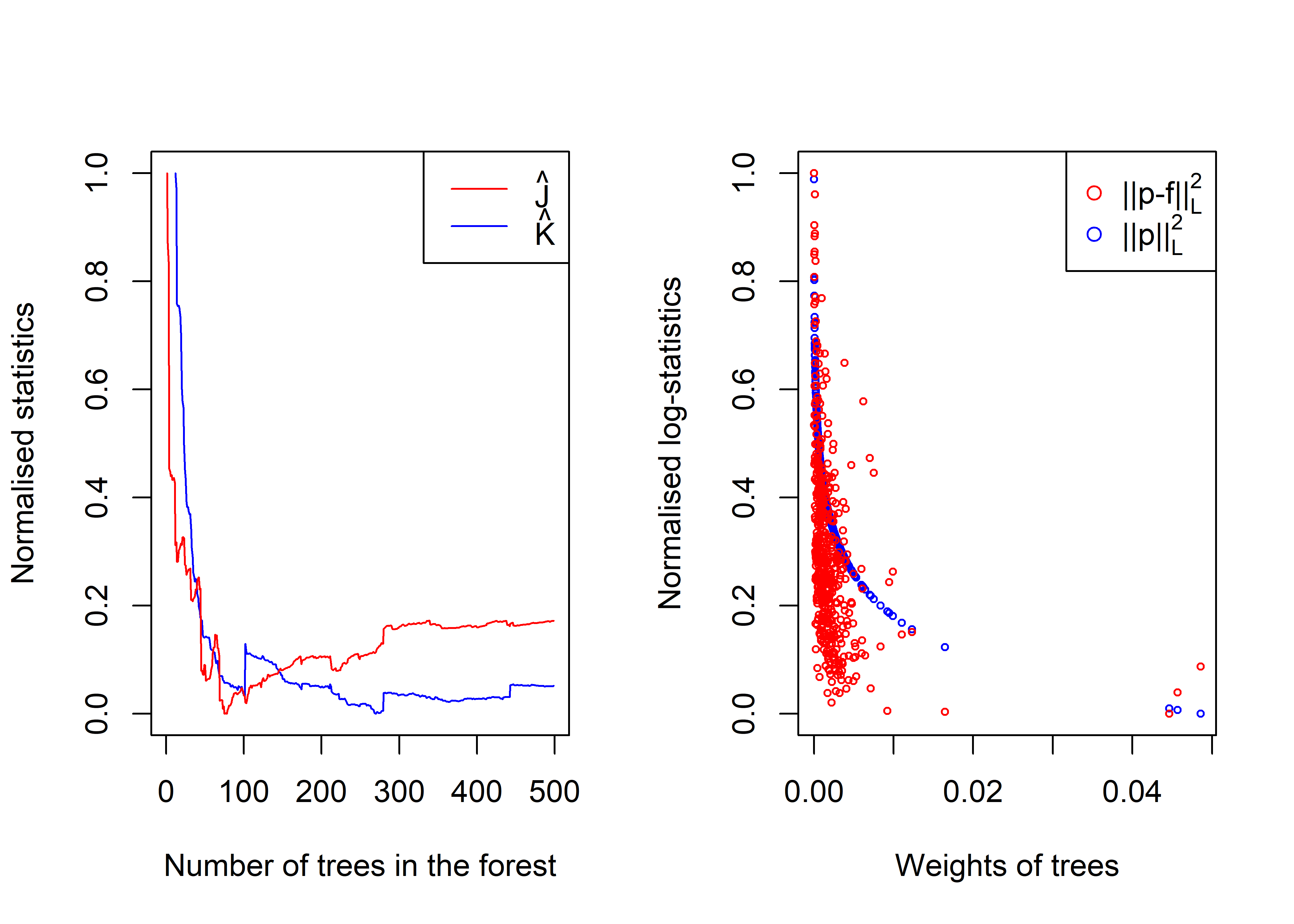} }}%
	\caption{(Dataset~\ref{data:clayton}) (a) Representation from the tree: in black, lower left, the input data. In red, upper-right, a simulation from the estimated tree. (b) Forest Statistics: on the left, $\hat{K}$ and $\hat{J}$ in function of the number of trees. On the right, the Integrated Constraint Influence and square norm of each tree against the weight of the tree in the forest.}%
	\label{fig:clayton_forest}%
\end{figure}

Dataset~\ref{data:clayton} is based on the Clayton copula, a commonly used dependence structure in many fields of application. The estimator developed in Algorithm~\ref{algo:CORT} has several options: the most important one is the inclusion, or not, of the localized dimension reduction through Algorithm~\ref{algo:dim_red}. Since here we have a completely independent dimension, this option is worth it: it reduces by a factor of $2$ the number of leaves, and hence the complexity of the model, by setting the same second edge $[0,1]$ to each leaf. Fig.~\ref{fig:clayton_forest} gives a representation of the tree and the statistics of the forest.  On the left of Fig.~\ref{fig:clayton_forest} (b), the convex, decreasing shape of the Kullback-Leibler divergence with respect to the number of trees shows that the generalization error of the forest decreases with the number of trees. The decreasing trend of the constraint influence and the square norm of trees with respect to the assigned weights by the forest, on the right of Fig.~\ref{fig:clayton_forest} (b) shows how the weighting procedure selected trees. Table~\ref{tab:clayton_table_dep_mes} shows dependence measures obtained from the different models. 

\begin{table}[H]
	\caption{\label{tab:clayton_table_dep_mes}Obtained dependence measures of several models on Dataset \ref{data:clayton}.  The first column is the goal, others are concurrent models.}
	\centering
	\begin{tabu} to \linewidth {>{\raggedright}X>{\raggedright}X>{\raggedright}X>{\raggedright}X>{\raggedright}X>{\raggedright}X>{\raggedright}X}
	\toprule
	  & Empirical & Cb(m=10) & Cb(m=5) & Beta & CORT & Bagged CORT\\
	\addlinespace[0.3em]
	\multicolumn{7}{l}{\Bold{Kendall Taus}}\\
	\hspace{1em}$\tau_{1,2}$ & -0.003 & 0.010 & 0.006 & -0.014 & 0.000 & -0.020\\
	\hspace{1em}$\tau_{1,3}$ & -0.796 & -0.750 & -0.673 & -0.799 & -0.780 & -0.493\\
	\hspace{1em}$\tau_{1,4}$ & 0.779 & 0.732 & 0.659 & 0.779 & 0.707 & 0.474\\
	\hspace{1em}$\tau_{2,3}$ & 0.015 & 0.010 & 0.011 & 0.029 & 0.000 & 0.031\\
	\hspace{1em}$\tau_{2,4}$ & -0.024 & -0.009 & -0.010 & -0.038 & 0.000 & -0.045\\
	\hspace{1em}$\tau_{3,4}$ & -0.775 & -0.728 & -0.654 & -0.773 & -0.695 & -0.566\\
	\addlinespace[0.3em]
	\multicolumn{7}{l}{\Bold{Spearman Rhos}}\\
	\hspace{1em}$\rho_{1,2}$ & -0.005 & 0.013 & 0.010 & -0.023 & 0.000 & -0.029\\
	\hspace{1em}$\rho_{1,3}$ & -0.934 & -0.915 & -0.868 & -0.936 & -0.926 & -0.648\\
	\hspace{1em}$\rho_{1,4}$ & 0.924 & 0.903 & 0.857 & 0.925 & 0.872 & 0.626\\
	\hspace{1em}$\rho_{2,3}$ & 0.023 & 0.014 & 0.016 & 0.045 & 0.000 & 0.047\\
	\hspace{1em}$\rho_{2,4}$ & -0.035 & -0.016 & -0.016 & -0.057 & 0.000 & -0.068\\
	\hspace{1em}$\rho_{3,4}$ & -0.922 & -0.901 & -0.853 & -0.922 & -0.862 & -0.735\\
	\bottomrule
	\end{tabu}
\end{table}

\begin{table}[H]
	\caption{\label{tab:clayton_table_bagging_result}Results of the bagging of each model on Dataset \ref{data:clayton}. Each row represents a different performance metric: in all cases, lower is better.}
	\centering
	\begin{tabu} to \linewidth {>{\raggedright}X>{\raggedright}X>{\raggedright}X>{\raggedright}X>{\raggedright}X>{\raggedright}X}
	\toprule
	  & Empirical & Cb(m=10) & Cb(m=5) & Beta & CORT\\
	$\hat{J}(\hat{c}_{\Gras{\omega}})$ & 0.00501 & -9.27 & -7.69 & 120 & \Bold{-54.6}\\
	$\hat{K}(\hat{c}_{\Gras{\omega}})$ & Inf & Inf & Inf & -0.582 & \Bold{-1.97}\\
	$\hat{M}(\hat{c}_{\Gras{\omega}})$ & 2.84e-05 & \Bold{2.51e-05} & 5.25e-05 & 3.38e-05 & 9.41e-05\\
	$\hat{N}(\hat{c}_{\Gras{\omega}})$ & -0.000666 & \Bold{-0.000669} & -0.000641 & -0.000657 & -0.000639\\
	\bottomrule
	\end{tabu}
\end{table}

The comparison of Table~\ref{tab:clayton_table_dep_mes} shows that, even if the CORT algorithm performs correctly, the forest tends to be biased in the dependence measures, toward more independence. On the other hand, the predictive performance of the model from Table~\ref{tab:clayton_table_bagging_result} is really high on the density-based estimates ($\hat{J}$ and $\hat{K}$), and is less good on the distribution function based versions ($\hat{M}$ and $\hat{N}$). 

The last dataset was produced based on a function $h_3$ defined by: $$h_3(\Gras{u}) = \left(u_1, \sin(2\pi u_1)-\frac{u_2}{\pi}, \left(1+\frac{u_3}{\pi^{2}}\right)\left(\frac{u_3}{2}\mathbb{1}_{\frac{1}{4} \ge u_1}  - \sin(\pi^{u_1})    \mathbb{1}_{\frac{1}{4} < u_1}\right)\right).$$

% % % % % % % % % % % % % % % % % % % % % % % % %  DATA FUNCDEP
% \noindent\begin{minipage}{0.45\textwidth}
% 	\begin{dataset}[Simulated functional]\label{data:funcdep}
% 		We chose to produce a voluntarily hard to estimate dependence structure, by applying the above function $h_3$ to uniformly drawn $3$-dimensional random vectors. The dataset is the ranks of $(h_3(\Gras{u}_i)))_{i \in 1,\ldots,500}$.
% 	\end{dataset}
% \end{minipage}
% \hfill%
% \begin{minipage}{0.45\textwidth}% adapt widths of minipages to your needs
% 	\centering\raisebox{\dimexpr \topskip-\height}{%
% 		\includegraphics[width=\textwidth]{funcdep_data.png}}
% 	\captionof{figure}{Pairs-plot of Dataset~\ref{data:funcdep}}
% \end{minipage}%
% % % % % % % % % % % % % % % % % % % % % % % % % 

% % % % % % % % % % % % % % % % % % % % % % % % DATA FUNCDEP
\begin{dataset}[Simulated functional]\label{data:funcdep}
	We chose to produce a voluntarily hard to estimate dependence structure, by applying $h_3$ to uniformly drawn $3$-dimensional random vectors. The dataset is the ranks of $(h_3(\Gras{u}_i)))_{i \in 1,\ldots,500}$.
\end{dataset}
% \begin{figure}[H]
% 	\centering
% 	\includegraphics[width=0.40\textwidth]{funcdep_data.png}
% 	\caption{Pairs-plot of Dataset~\ref{data:funcdep}.}
% \end{figure}
% % % % % % % % % % % % % % % % % % % % % % % % 

\begin{figure}[H]
	\centering
	\subfloat[]{{\includegraphics[width=0.40\textwidth]{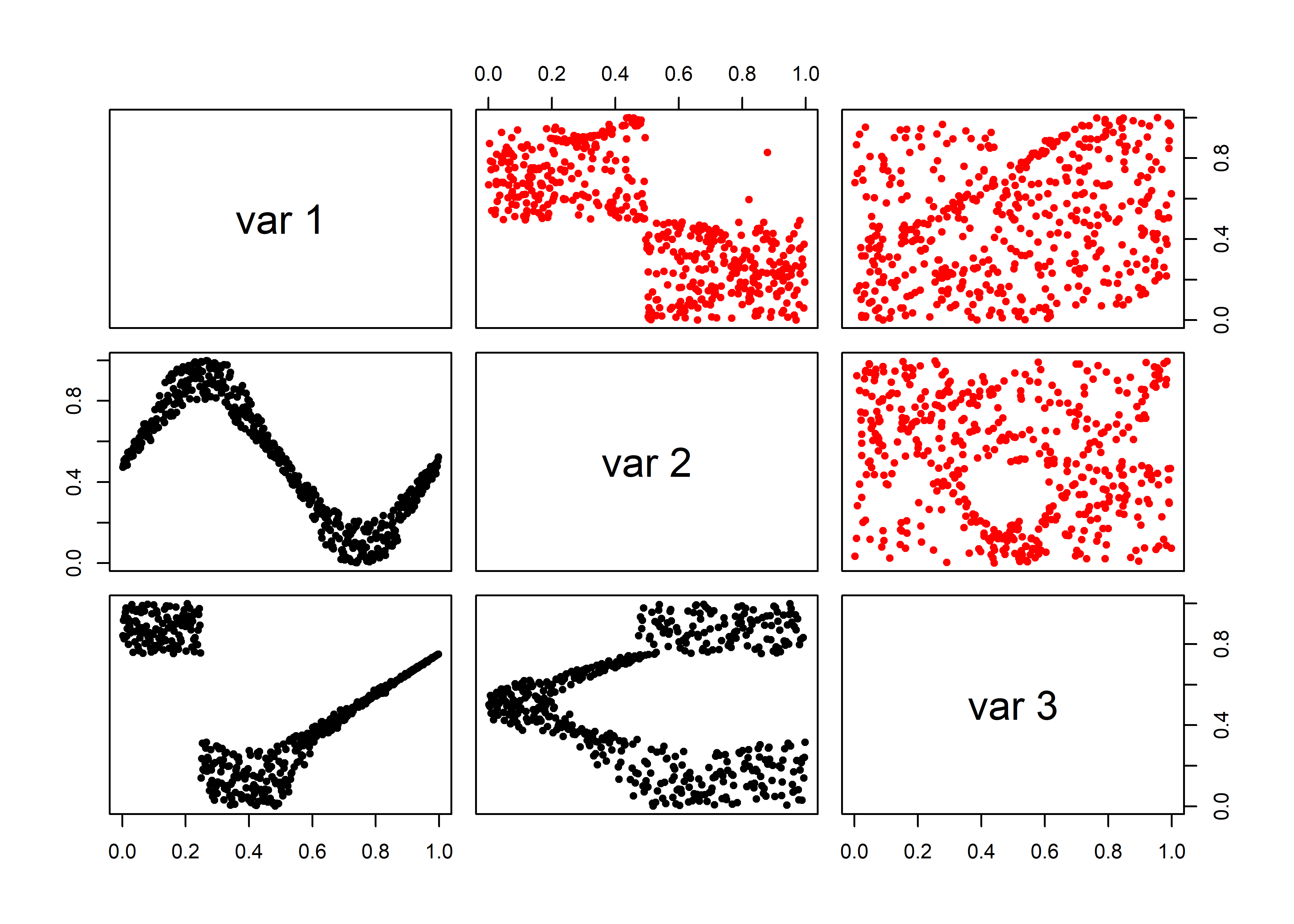} }}%
	\qquad
	\subfloat[]{{\includegraphics[width=0.50\textwidth]{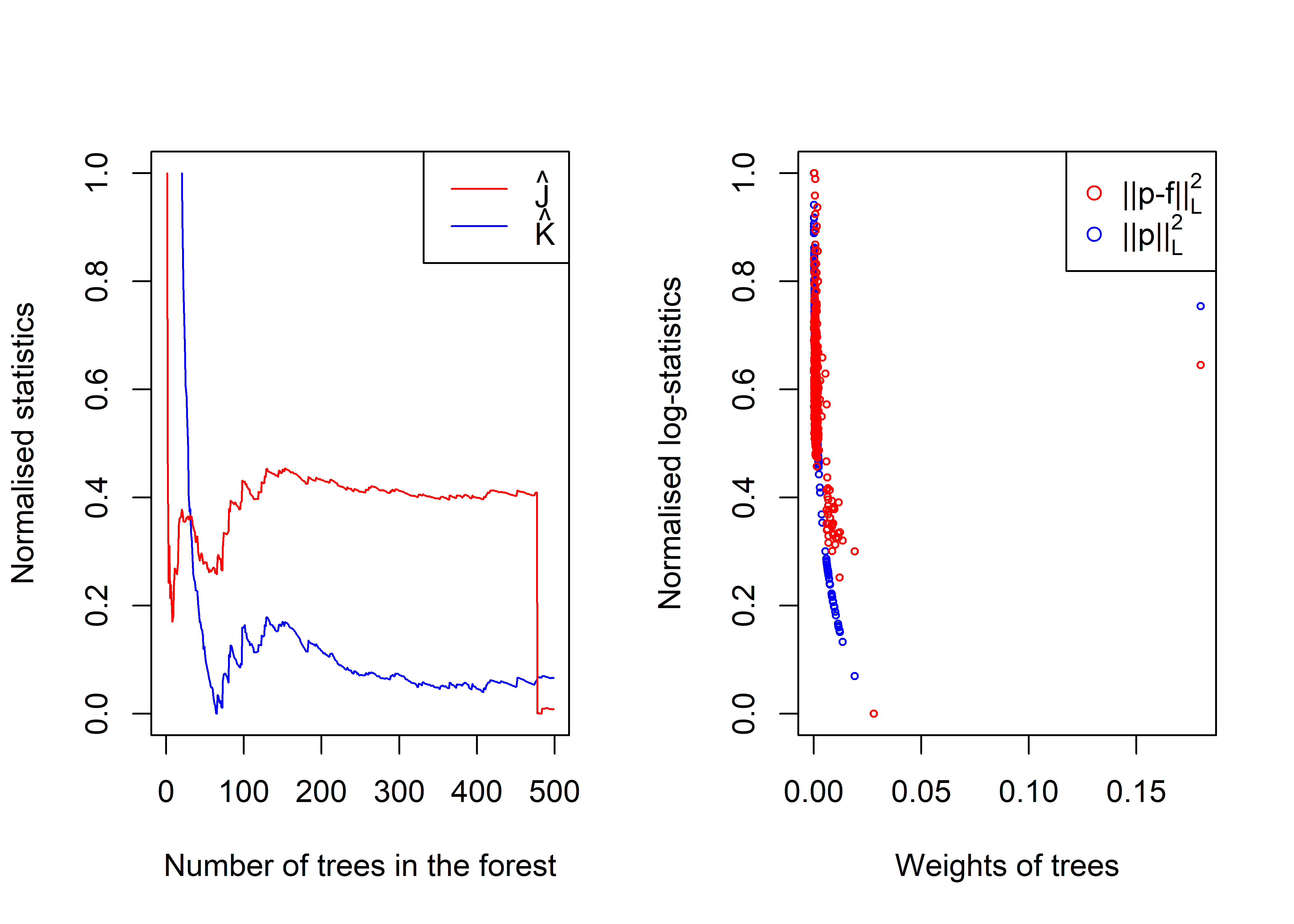} }}%
	\caption{(Dataset~\ref{data:funcdep}) (a) Representation from the tree: in black, lower left, the input data. In red, upper-right, a simulation from the estimated tree. (b) Forest Statistics: On the left, $\hat{K}$ and $\hat{J}$ in function of the number of trees. On the right, the Integrated Constraint Influence and square norm of each tree against the weight of the tree in the forest.}%
	\label{fig:funcdep_forest}%
\end{figure}

Fig.~\ref{fig:funcdep_forest} shows the CORT estimator and the bagging statistics on this dataset. 

Although the estimation on the second and third marginals is not very good, the statistics of the forests are good, showing that the estimator get better and better while adding trees. Indeed, on Fig.~\ref{fig:funcdep_forest} (b), left, we observe decreasing out-of-bag errors, but we also observe on Fig.~\ref{fig:funcdep_forest} (b), right, a very skewed and fat-tailed density for the constraint influence, meaning that certain trees did not produce very good partitions.

\begin{table}[H]
	\caption{\label{tab:funcdep_table_dep_mes}Obtained dependence measures of several models on Dataset \ref{data:funcdep}. The first column is the goal, others are concurrent models.}
	\centering
	\begin{tabu} to \linewidth {>{\raggedright}X>{\raggedright}X>{\raggedright}X>{\raggedright}X>{\raggedright}X>{\raggedright}X>{\raggedright}X}
	\toprule
	  & Empirical & Cb(m=10) & Cb(m=5) & Beta & CORT & Bagged CORT\\
	\addlinespace[0.3em]
	\multicolumn{7}{l}{\Bold{Kendall Taus}}\\
	\hspace{1em}$\tau_{1,2}$ & -0.495 & -0.491 & -0.471 & -0.500 & -0.458 & -0.451\\
	\hspace{1em}$\tau_{1,3}$ & 0.089 & 0.042 & -0.002 & 0.084 & 0.171 & 0.048\\
	\hspace{1em}$\tau_{2,3}$ & 0.005 & 0.015 & 0.004 & 0.007 & -0.109 & 0.013\\
	\addlinespace[0.3em]
	\multicolumn{7}{l}{\Bold{Spearman Rhos}}\\
	\hspace{1em}$\rho_{1,2}$ & -0.743 & -0.737 & -0.715 & -0.747 & -0.713 & -0.674\\
	\hspace{1em}$\rho_{1,3}$ & -0.156 & -0.154 & -0.146 & -0.159 & 0.246 & 0.060\\
	\hspace{1em}$\rho_{2,3}$ & -0.019 & 0.012 & -0.001 & -0.010 & -0.177 & 0.020\\
	\bottomrule
	\end{tabu}
\end{table}

Table~\ref{tab:funcdep_table_dep_mes} shows that the bivariate projections were not all treated as well as others by the CORT algorithm: the values of $\tau_{1,2}$ and $\rho_{1,2}$ are surprisingly quite good compared to $\tau_{1,3},\tau_{2,3},\rho_{1,3},\rho_{2,3}$, for the CORT estimator. Hopefully, the bagging corrects this bias quite correctly. 

\begin{table}[H]
	\caption{\label{tab:funcdep_table_bagging_result}Results of the bagging of each model on Dataset \ref{data:funcdep}. Each row represents a different performance metric: in all cases, lower is better.}
	\centering
	\begin{tabu} to \linewidth {>{\raggedright}X>{\raggedright}X>{\raggedright}X>{\raggedright}X>{\raggedright}X>{\raggedright}X}
	\toprule
	  & Empirical & Cb(m=10) & Cb(m=5) & Beta & CORT\\
	$\hat{J}(\hat{c}_{\Gras{\omega}})$ & 0.002 & -15.7 & -7.22 & -23 & \Bold{-23.3}\\
	$\hat{K}(\hat{c}_{\Gras{\omega}})$ & Inf & -2.52 & -1.79 & \Bold{-3.12} & -1.72\\
	$\hat{M}(\hat{c}_{\Gras{\omega}})$ & \Bold{1.69e-05} & 0.000172 & 0.000736 & 4.76e-05 & 0.00233\\
	$\hat{N}(\hat{c}_{\Gras{\omega}})$ & \Bold{-0.0201} & -0.0199 & -0.0194 & -0.0201 & -0.0177\\
	\bottomrule
	\end{tabu}
\end{table}

Finally, the predictive performance from Table~\ref{tab:funcdep_table_bagging_result} is still two-sided: the density-based results are quite good, but the distribution-function based ones are not very good. 
% The box plot of Fig.~\ref{fig:funcdep_Box plot} confirms this analysis.

% \begin{figure}[H]
% 	\centering
% 	\subfloat[Box plot of log-normalized $\hat{P}$]{{\includegraphics[width=0.40\textwidth]{funcdep_boxplotP_lognorm.png}}}%
% 	\qquad
% 	\subfloat[Box plot of log-normalized $\hat{Q}$]{{\includegraphics[width=0.40\textwidth]{funcdep_boxplotQ_lognorm.png}}}%
% 	\caption{(Dataset~\ref{data:funcdep})  Box plots of resulting errors for $20$ resamples for each model (the lower the better). $\hat{P}$ is focused on the distribution function and $\hat{Q}$ on the density.}%
% 	\label{fig:funcdep_Box plot}%
% \end{figure}

More details and experiments are available in the supplementary section. 

\section{Conclusion\label{sec:6}}

From a simple density estimation procedure designed by~\cite{ram2011}, we constructed a piecewise constant, tree-shaped, recursive copula density estimator. We computed several closed-form expressions for this estimator, and we gave an asymptotic result.

If, intuitively, constraining the space of potential weighting solutions will help the convergence of optimization routines, the copula constraints forced us to split the space in more than one dimension, making the resulting estimation procedure complex with the increasing dimension. The localized dimension reduction procedure helps to reduce the complexity.

The CORT estimator has good generalization performance and is straightforward to use since it does not have restrictive hypothesis on the true dependence structure. Although the implementation we provide is very fast, a balance between computation time and precision is available in the number of trees used in the bagging procedure. However, more work needs to be done to correct defaults of the splitting procedure, which is not able to understand certain kinds of dependence structures. 

\section*{Acknowledgments}

We are particularly grateful to the two referees, the associate editor and the editor who all made useful comments and suggestions, improving the manuscript. We would like to thank SCOR SE for funding this work through a CIFRE grant. Any errors are ours.

\section*{Supplementary material}

On each dataset from Section~\ref{sec:5}, we ran an additional experiment. To observe the predictive performance of each model, we designed a cross-validation procedure: on $20$ resamples of each dataset, we computed the Cramer-Von-Mises and ISE errors on test samples, given respectively by:
\begin{align*}
\hat{P}(\hat{D}) &= \sum_{i \in T} \left(\hat{D}(\Gras{u}_i) - \hat{C}(\Gras{u}_i)\right)^2,\qquad \hat{Q}(\hat{d}) = \norme{\hat{d}} - \frac{2}{\lvert T \rvert}\sum_{i \in T} \hat{d}(\Gras{u}_i),
\end{align*}
where $T$ is the test dataset used to obtain $\hat{D}$, a given copula estimator with density $\hat{d}$.  Note that $\hat{P}$ is focused on the distribution function and $\hat{Q}$ on the density. Below are the resulting box plots for each of the datasets. 

\begin{figure}[H]
	\centering
	\subfloat[Box plot of log-normalized $\hat{P}$]{{\includegraphics[width=0.40\textwidth]{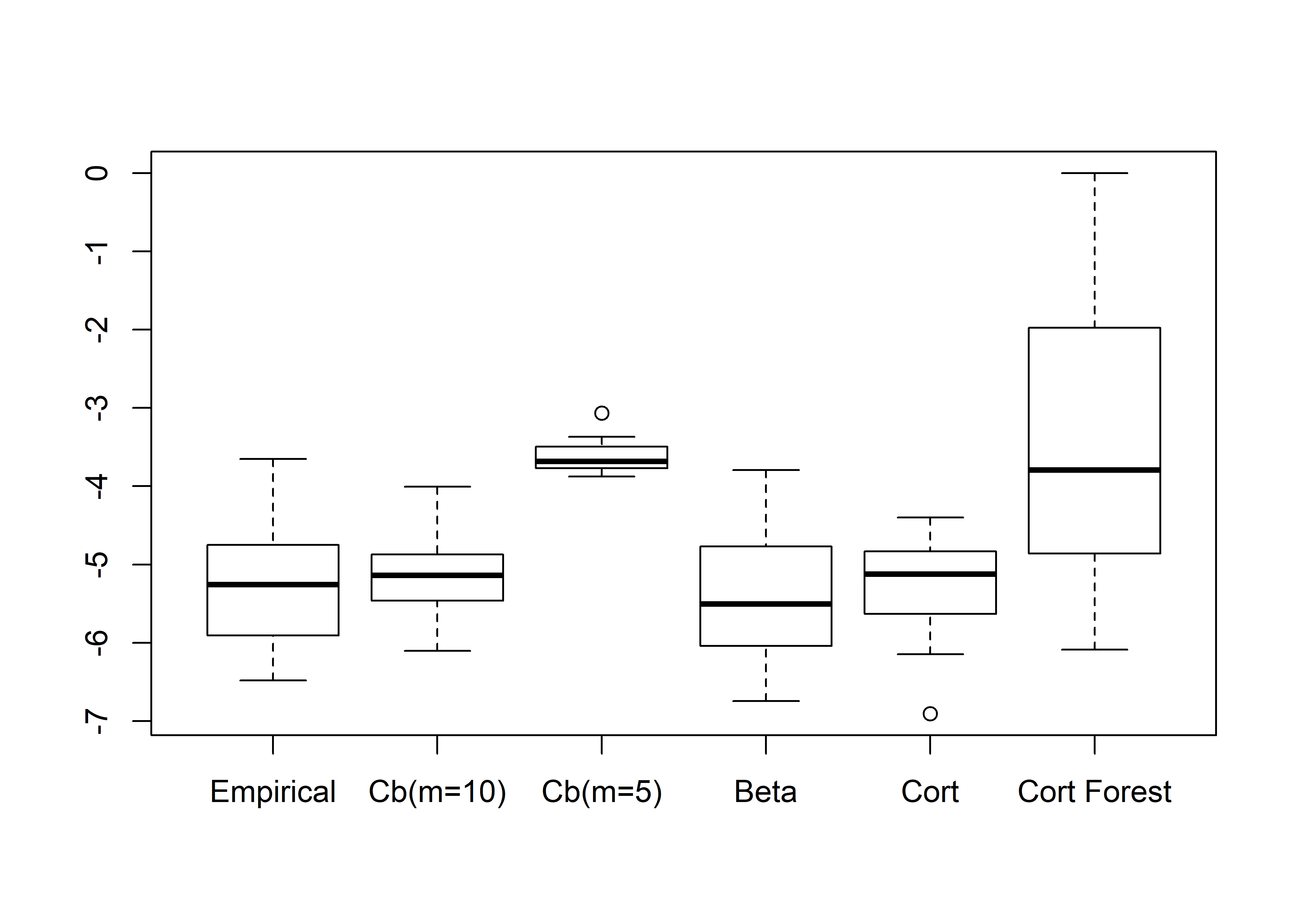}}}%
	\qquad
	\subfloat[Box plot of log-normalized $\hat{Q}$]{{\includegraphics[width=0.40\textwidth]{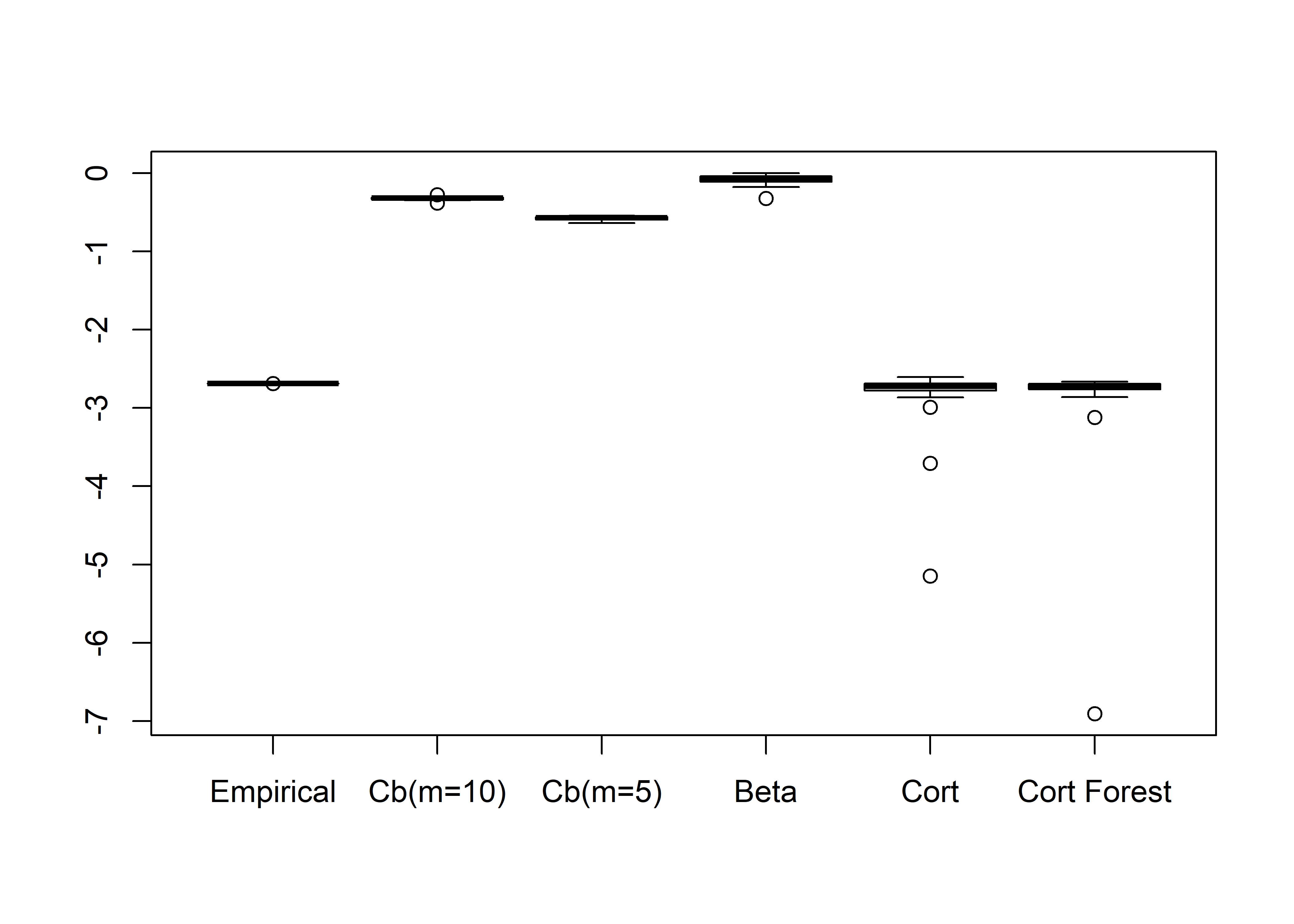}}}%
	\caption{(Dataset~\ref{data:recoveryourself})  Box plots of resulting errors for $20$ resamples for each model (the lower the better). $\hat{P}$ is focused on the distribution function and $\hat{Q}$ on the density.}%
	\label{fig:recoveryourself_Box plot}%
\end{figure}

Fig.~\ref{fig:recoveryourself_Box plot} gives a box plot of $\hat{P}$ Cramer-Von-Mises errors on the first dataset.

On Fig.~\ref{fig:recoveryourself_Box plot}, note that smaller values of $\hat{P}$ (d.f. based) and $\hat{Q}$ (density based) mean a better model. We observe that, although the bagging procedure is not worth it, the CORT estimator is very good on this example, both for density and d.f. estimation. Unfortunately, it is not always the case, as shown by the experiments we did on Dataset~\ref{data:impossible}.

\begin{figure}[H]
	\centering
	\subfloat[Box plot of log-normalized $\hat{P}$]{{\includegraphics[width=0.40\textwidth]{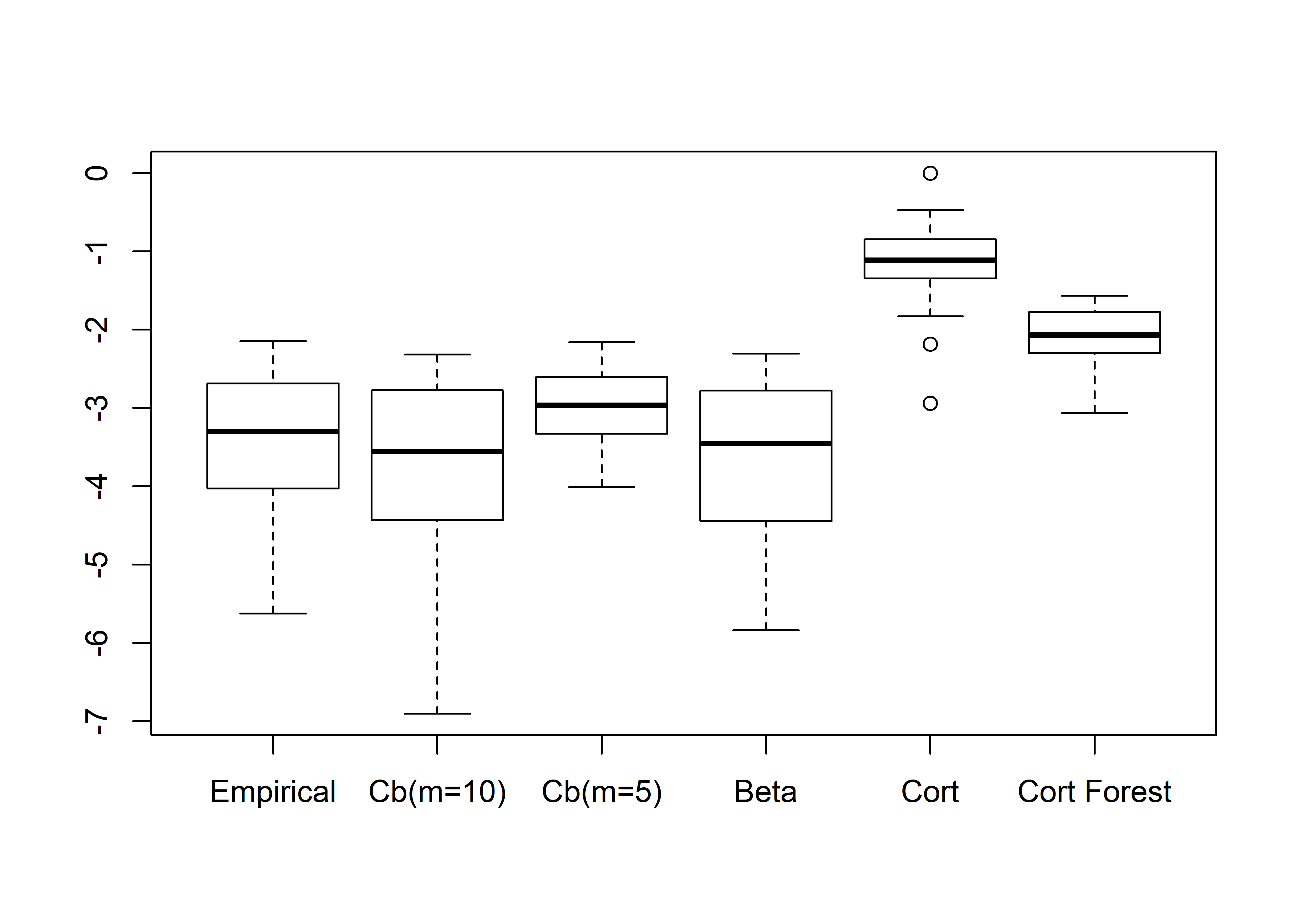}}}%
	\qquad
	\subfloat[Box plot of log-normalized $\hat{Q}$]{{\includegraphics[width=0.40\textwidth]{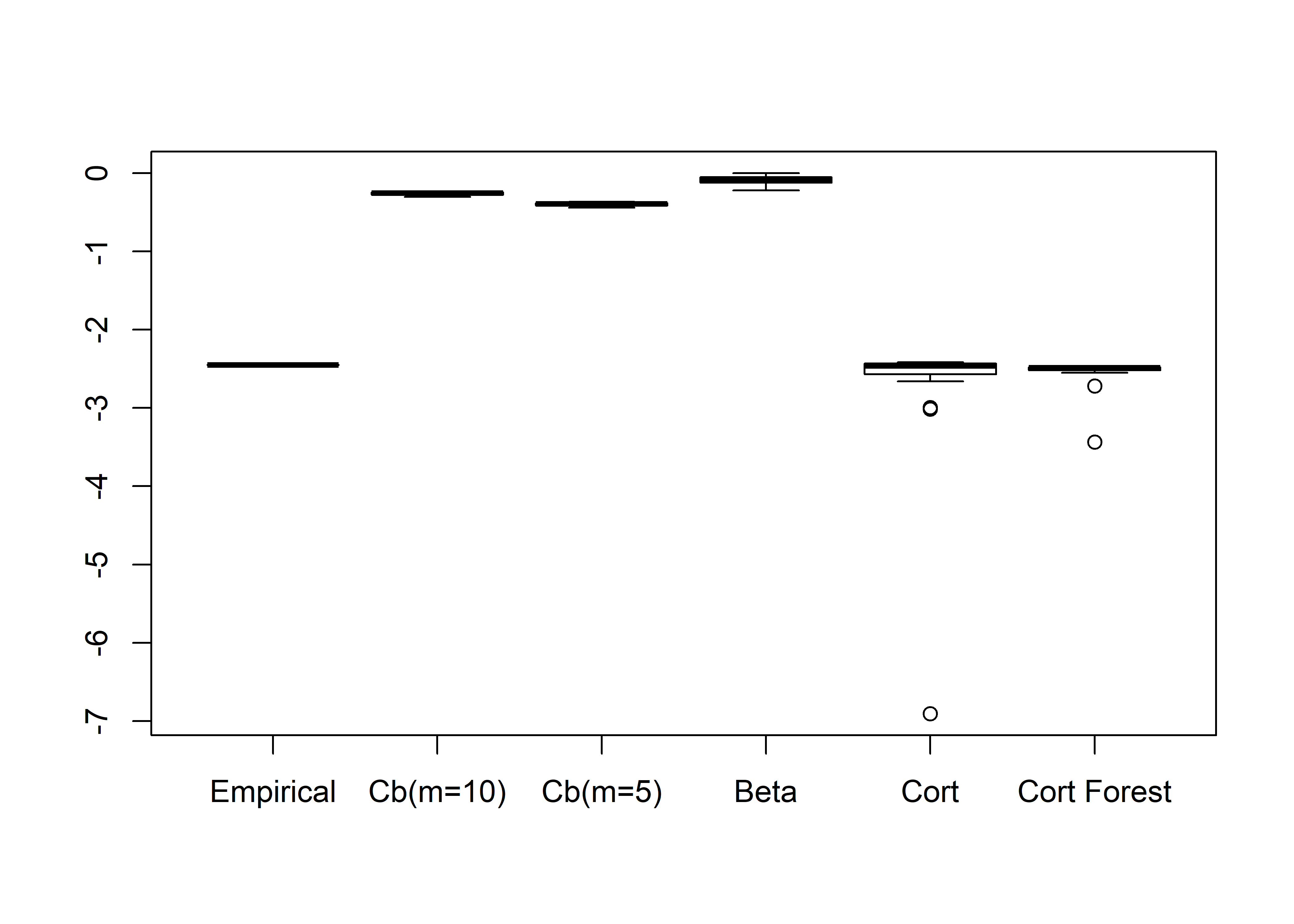}}}%
	\caption{(Dataset~\ref{data:impossible})  Box plots of resulting errors for $20$ resamples for each model (the lower the better). $\hat{P}$ is focused on the distribution function and $\hat{Q}$ on the density.}%
	\label{fig:impossible_Box plot}%
\end{figure}
The box plot of $\hat{P}$ error on Dataset~\ref{data:impossible} confirms the analysis we already had: in this case, the empirical beta copula performs a lot better.

For Dataset~\ref{data:clayton}, however, the box plot of $\hat{P}$ and $\hat{Q}$ in Fig.~\ref{fig:clayton_Box plot} confirms the performance of the estimator:

\begin{figure}[H]
	\centering
	\subfloat[Box plot of log-normalized $\hat{P}$]{{\includegraphics[width=0.40\textwidth]{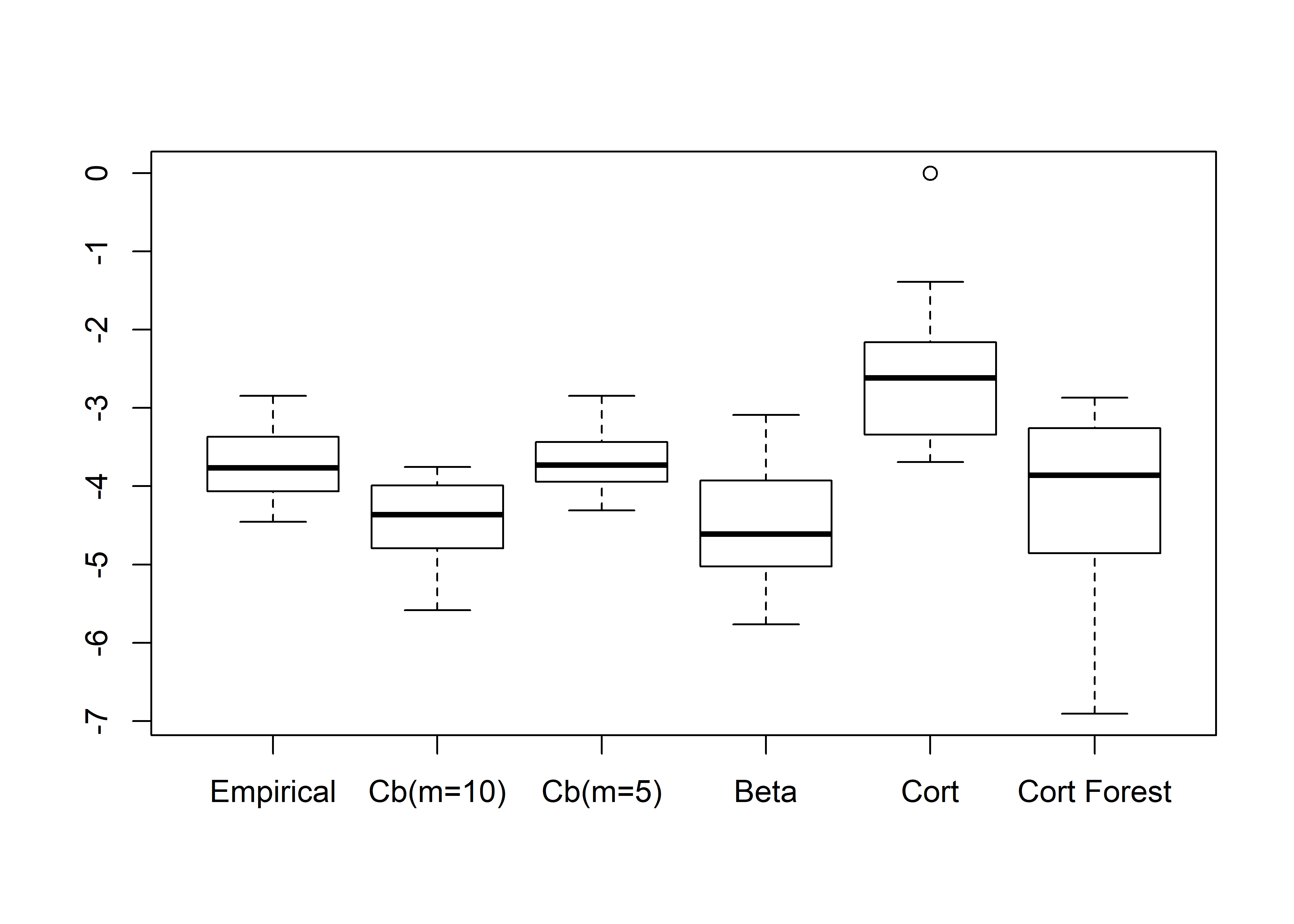}}}%
	\qquad
	\subfloat[Box plot of log-normalized $\hat{Q}$]{{\includegraphics[width=0.40\textwidth]{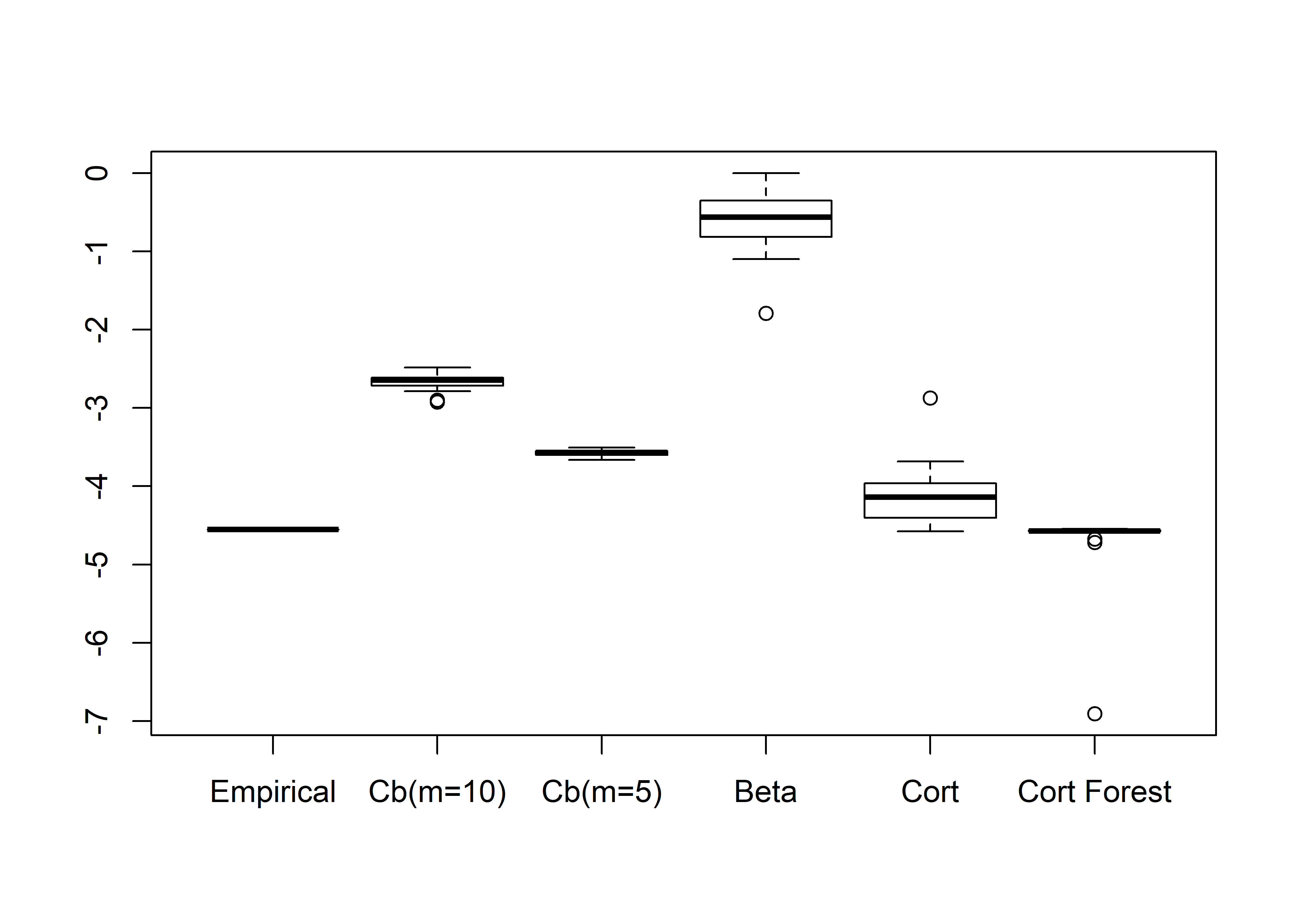}}}%
	\caption{(Dataset~\ref{data:clayton})  Box plots of resulting errors for $20$ resamples for each model (the lower the better). $\hat{P}$ is focused on the distribution function and $\hat{Q}$ on the density.}%
	\label{fig:clayton_Box plot}%
\end{figure}

\begin{figure}[H]
	\centering
	\includegraphics[width=0.7\textwidth]{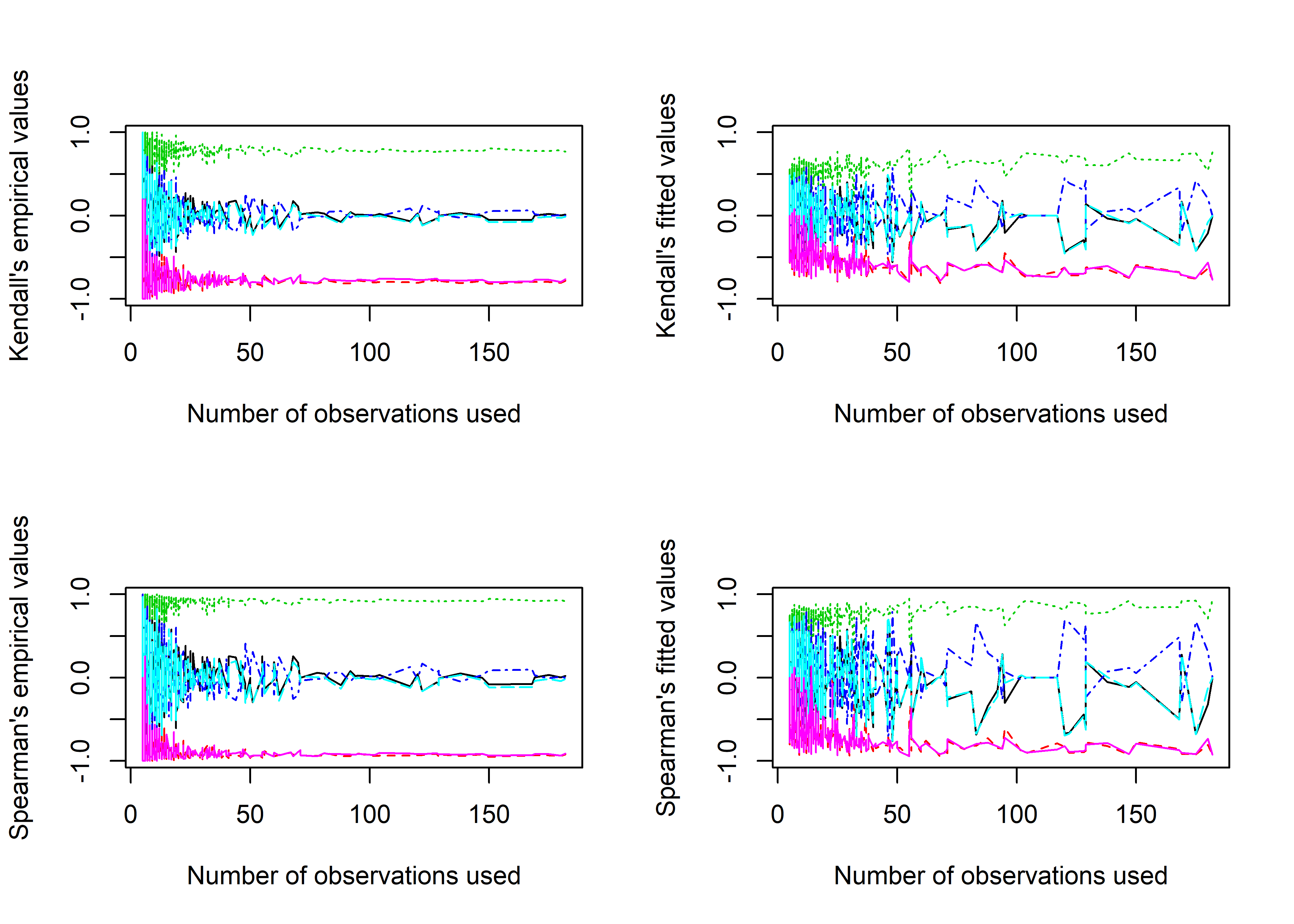}
	\caption{(Dataset~\ref{data:clayton}) Burn in of dependence measures: on the top row, Kendall's tau; on the bottom row, Spearman's rho; on the left column, empirical values from subsamples of increasing sizes; on the right column, values obtained by the fitted CORT estimators on the same subsamples. The size of the subsamples is in abscissa. Each line type corresponds to a couple of variables.}
	\label{fig:clayton_values_kendall_spearmann}
\end{figure}

To understand more precisely what happened with Kendall's tau and Spearman's rho on Dataset~\ref{data:clayton}, we ran a burn-in experiment: we fitted trees on subsamples of increasing size. We can then observe the burn in Kendall taus and Spearman rhos, represented on Figure~\ref{fig:clayton_values_kendall_spearmann}.  We see that the fitted values, on the right of Fig.~\ref{fig:clayton_values_kendall_spearmann}, are convergent but biased compared to the empirical observed values of the dependence measures directly computed on resamples of the dataset, on the left. We also observe the high variance of the estimator, which is one of the good reasons to use a bagging procedure.

\begin{figure}[H]
	\centering
	\subfloat[Box plot of log-normalized $\hat{P}$]{{\includegraphics[width=0.40\textwidth]{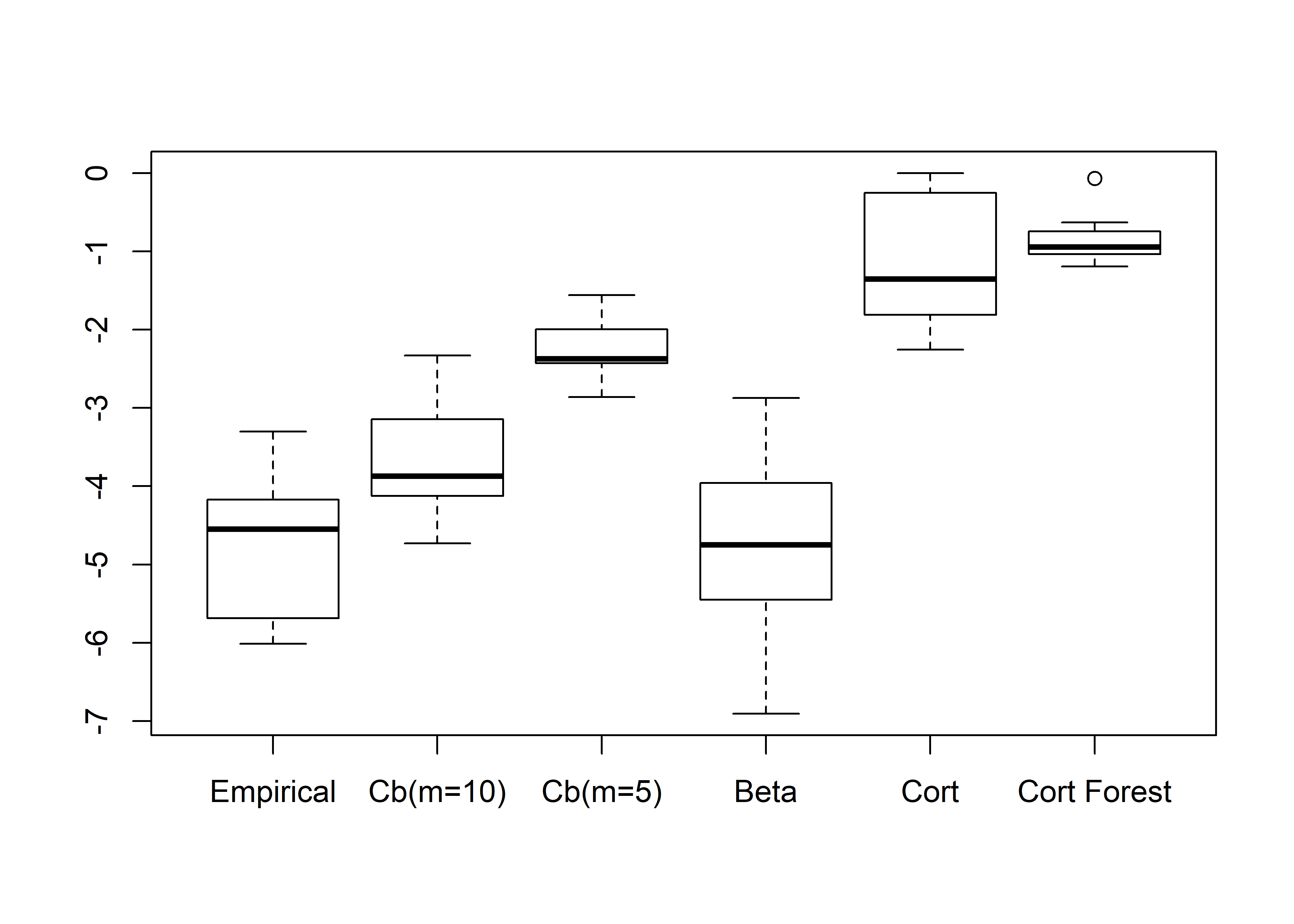}}}%
	\qquad
	\subfloat[Box plot of log-normalized $\hat{Q}$]{{\includegraphics[width=0.40\textwidth]{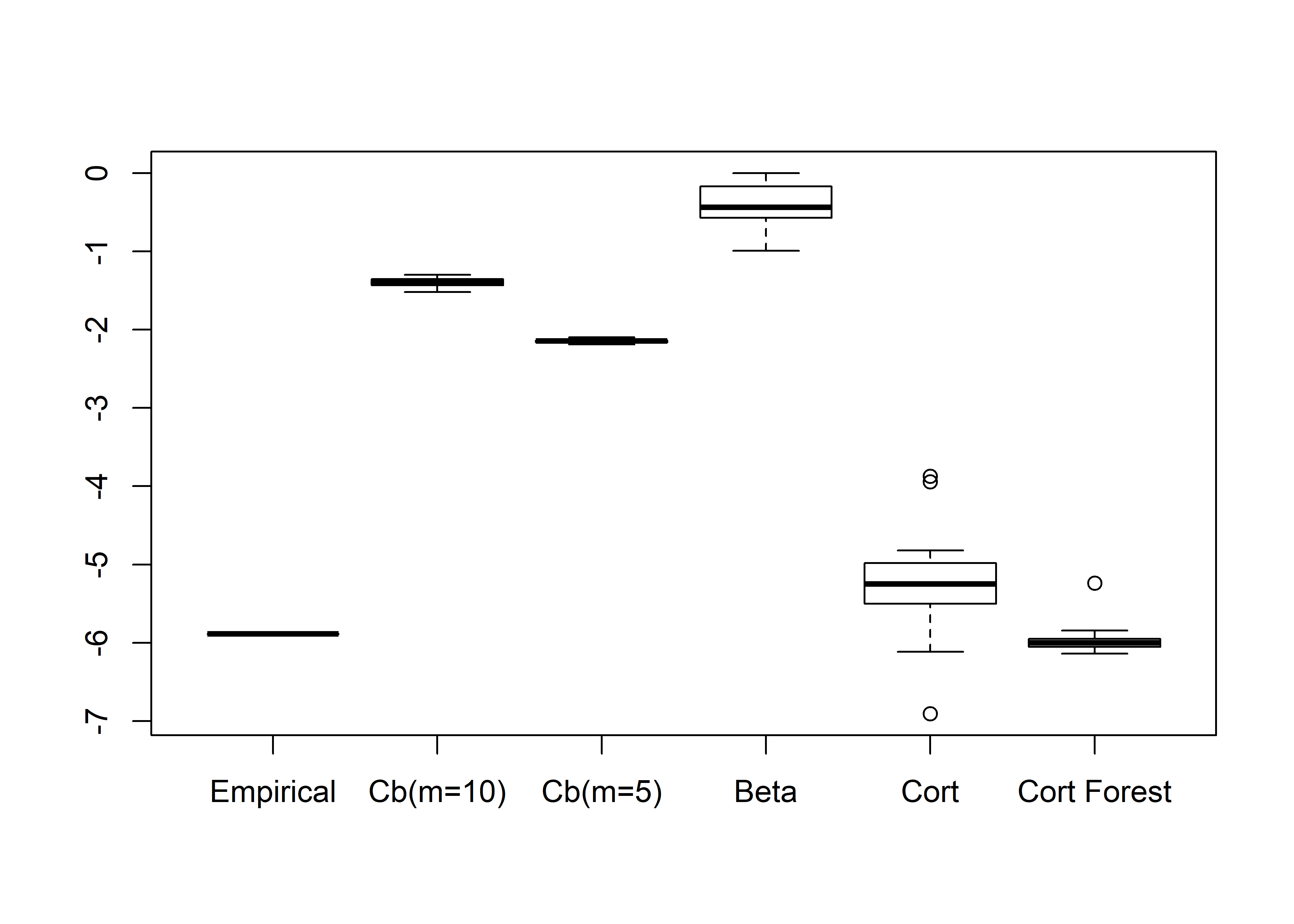}}}%
	\caption{(Dataset~\ref{data:funcdep})  Box plots of resulting errors for $20$ resamples for each model (the lower the better). $\hat{P}$ is focused on the distribution function and $\hat{Q}$ on the density.}%
	\label{fig:funcdep_Box plot}%
\end{figure}

Finally, the box plots on Dataset~\ref{data:funcdep} are given in Fig.~\ref{fig:funcdep_Box plot}.

\begin{proof} of Proposition~\ref{ppt:comon_dep_mes}. We have: 
	
	\begin{align*}
	\int c_{\Gras{p},\L}(\Gras{u}) \, C_{\Gras{p},\L}(\Gras{u}) \;d\Gras{u} &= \sum\limits_{\ell \in \L}\sum\limits_{k \in \L}\frac{p_\ell p_k}{\lambda(\ell) \lambda(k)}\int\limits_{\ell} \lambda([0,\Gras{u}]\cap k)\;d\Gras{u} \\
	&= \sum\limits_{\ell \in \L}\sum\limits_{k \in \L}\frac{p_\ell p_k}{\lambda(\ell) \lambda(k)}\prod\limits_{i=1}^{d} \int\limits_{\ell_i} \lambda([0,u_i]\cap k_i)\;du_i
	\end{align*}
	
	where $\ell_i$ denotes the projection of $\ell$ onto the dimension $i$. Denote now $g(u,c,d) = \lambda([0,u]\cap [c,d])$, and remark that: 
	
	$$g(u,c,d) = \lambda([0,u]\cap [c,d]) = \begin{cases} 0, & u \le c\\u-c, & c< u < d\\d-c,& d \le u\end{cases}$$
	
	Denote furthermore $G(u,a,b,c,d) = \int\limits_{a}^{b} g(u,c,d) du$. We have: 
	
	\begin{align*}
	G(u,a,b,c,d) &= \int_{a}^{b}0\,\mathbb{1}_{u<c}\,du &&+\int_{a}^{b}\left(u-c\right)\mathbb{1}_{c<u<d}\,du &&+\int_{a}^{b}\left(d-c\right)\mathbb{1}_{d<u}\,du \\
	&= &&\int_{a\vee c}^{b \wedge d}\left(u-c\right)\mathbb{1}_{c<u<d}\,du &&+\int_{a\vee d}^{b}\left(d-c\right)\mathbb{1}_{d<u}\,du\\
	&= &&\left(\frac{(b\wedge d)^2}{2} - (b\wedge d)c\right) - \left(\frac{(a\vee c)^2}{2} - (a\vee c)c\right) &&+(d-c)\left(b - (a\vee d)\right),\\
	\end{align*}

	where the two integrals are simply integrals of linear functions, regarding the univariate slack variable $u$. Hence, $$G(u,a,b,c,d)  = \frac{(b\wedge d)^2 - (a\vee c)^2}{2} + c\left((a\vee c) - (b\wedge d)\right)+(d-c)\left(b - (a\vee d)\right)\text{,}$$ which provides the wanted expression for $\tau$. For $\rho$, one only needs to compute the integral expression:
	
	\begin{align*}
	\int C_{\Gras{p},\L}(u)\, du &= \sum\limits_{\ell \in \L} \frac{p_\ell}{\lambda(\ell)} \int \lambda([0,\Gras{u}]\cap \ell) \,du\\
	&= \sum\limits_{\ell \in \L} \frac{p_\ell}{\lambda(\ell)} \prod\limits_{i=1}^{d}  \int\limits_{0}^{a_i} 0 \,du_i +  \int\limits_{a_i}^{b_i} \left(u_i-a_i\right) \,du_i + \int\limits_{b_i}^{1} \left(b_i-a_i\right) \,du_i \\
	& = \sum\limits_{\ell \in \L} \frac{p_\ell}{\lambda(\ell)} \prod\limits_{i=1}^{d}  0 +  \frac{b_i^2 - a_i^2}{2} - (b_i - a_i)a_i + (1-b_i)(b_i - a_i) \\
	& = \sum\limits_{\ell \in \L} \frac{p_\ell}{\lambda(\ell)} \prod\limits_{i=1}^{d} \frac{1}{2}(b_i - a_i)(2 - b_i - a_i)
	\end{align*}
	
	which concludes the argument since $\lambda(\ell) = \prod\limits_{i=1}^{d}(b_i - a_i)$.
	
\end{proof}

\bibliographystyle{myjmva}
\bibliography{everything}
%\begin{thebibliography}
%\end{thebibliography}

\end{document}